\documentclass[journal,article,accept,moreauthors,12pt,a4paper]{mdpi} 

\lastpage{4699}
\doinum{10.3390/e15114668}
\pubvolume{15}
\pubyear{2013}
\history{Received: 3 August 2013; in revised form: 11 September 2013 / Accepted: 17 October 2013 / \\
Published: 31 October 2013}

\usepackage{soul}

\usepackage{graphicx}
\usepackage{amsmath} 
\usepackage{amssymb} 
\usepackage{mathrsfs} 
\newcommand{\E}{{\mathbb{E}}}

\newcommand{\m}{{|Z|}}
\newcommand{\n}{{\vec{n}}}
\newcommand{\kk}{{\vec{k}}}
\newcommand{\D}{{\mathscr{D}}}
\newcommand{\SSS}{{\mathscr{S}}}

\usepackage{color} 
\RequirePackage[dvipsnames,usenames]{xcolor}

\newcommand{\qed}[0]{{$\blacksquare$}}
\newtheorem{theorem}{Theorem}

\newtheorem{corollary}[theorem]{Corollary}

\newtheorem{definition}{Definition}

\newtheorem{lemma}[theorem]{Lemma}

\newtheorem{proposition}[theorem]{Proposition}

\Title{Estimating Functions of Distributions Defined over Spaces of Unknown Size} 

\Author{David H. Wolpert $^{1,}$* and Simon DeDeo $^{1,2}$}

\address{%
$^{1}$ Santa Fe Institute, 1399 Hyde Park Rd., Santa Fe, NM 87501, USA; E-Mail: simon@santafe.edu\\
$^{2}$ School of Informatics and Computing, Indiana University, 901 E 10th St, Bloomington, IN~47408,~USA \\
}

\corres{E-Mail: david.h.wolpert@gmail.com.}

\abstract{We consider Bayesian estimation of information-theoretic quantities
from data, using a Dirichlet prior. Acknowledging the uncertainty of the event space size $m$
and the Dirichlet prior's concentration parameter $c$, we treat both 
as random variables set by a hyperprior. We show that the associated hyperprior, $P(c, m)$, obeys a simple ``Irrelevance of 
Unseen Variables'' (IUV) desideratum iff $P(c, m) = P(c) P(m)$. 
Thus, requiring IUV greatly reduces the number of degrees of freedom of the hyperprior.
Some information-theoretic quantities can be expressed multiple ways,
in terms of different event spaces, e.g., mutual information. With all hyperpriors (implicitly) used in earlier work,
different choices of this event space lead to different posterior expected values
of these information-theoretic quantities. We show that there is no such dependence on
the choice of event space for a hyperprior that obeys IUV.
We also derive a result that allows us to exploit IUV to greatly simplify calculations, like
the posterior expected mutual information or posterior expected multi-information.
We also use computer experiments to favorably compare an IUV-based estimator of entropy
to three alternative methods in common use.
We end by discussing how seemingly innocuous changes to the formalization of an 
estimation problem can substantially
affect the resultant estimates of posterior expectations.
}


\begin{document}

%
%
%
%
%
%
%
%
%

\newpage


\section{Background}
\label{sec:background}

A central problem of statistics is estimating a functional of
a probability distribution $\rho$ from a dataset, $\n$, of $N$ independent, identically distributed (IID
)
samples of $\rho$. A simple example of this problem is estimating the
mean $Q(\rho)$ of a distribution, $\rho$, from a set of IID samples of that distribution.
In this example, the functional $Q(.)$ depends linearly on $\rho$.
More challenging versions of this problem arise when the functional $Q(.)$
is nonlinear. In particular, recent decades have seen a lot of work
on estimating information-theoretic functionals~\cite{coth91,mack03} of a distribution
from a set of IID samples of that distribution. Examples include estimating 
the Shannon entropy of the distribution, its mutual information, \textit{etc}., from such a set of 
samples and, in particular, using the bootstrap to estimate
associated error bars~\cite{paninski2003estimation,grassberger2003entropy,korber1993covariation}.

This work has concentrated on the case where the
event space being sampled is countable.
In addition, much of it has used non-Bayesian approaches.
The first work that addressed the problem using Bayesian techniques was
the sequence of papers~\cite{wowo95a,wowo95b,wowo_erratum} (hereafter abbreviated
as WW
), followed by similar work in~\cite{hutter2002distribution} (see Appendix A 
for a list of corrections to some algebraic mistakes in~\cite{wowo95a};
a careful analysis of the numerical implementation of the formulas in WW can be found in~\cite{huka13}).
This work showed how to calculate posterior moments of several
nonlinear functionals of the distribution. Such moments provide both
the Bayes-optimal estimate of the functional (assuming a quadratic loss
function) and an error bar in that estimate. In particular, WW provided closed-form expression for 
the posterior first and
second moments of entropy, mutual information and Kullback-Leibler
distance, in addition to various \linebreak non-information-theoretic quantities, like
covariance.

Write the space of possible events as $Z$ with elements written as $z$
and distributions over $Z$ as $\rho$.
For tractability reasons, WW used a Dirichlet prior over the associated simplex of 
possible distributions,
$P(\rho) \propto \prod_z \rho(z)^{c L(z) - 1}$.
(In the literature, the constant, $c$,
is sometimes called a {concentration parameter}, and $L$ is sometimes called a 
{baseline distribution}.) 

In WW, $Z$ was taken to be fixed (not a random variable), $L(z)$ was taken to be uniform over $Z$
and $c$ was taken to equal $|Z|$, the size of $Z$.
This choice of a Dirichlet prior over $\rho$ with uniform $L$ has been the basis of all
subsequent work on Bayesian estimates
of information-theoretic functionals of distributions~\cite{Nemenman:2002fk,archer2012Bayesian}.
(Note though that recently, there has been some investigation of the extension
of this work to mixture-of-Dirichlet distributions~\cite{Nemenman:2002fk} and Dirichlet / Pittman-Yor 
processes~\cite{archer2012Bayesian}).)
However, there has not been such consensus concerning $c$. 

Although the choice of $c = |Z|$ was not explicitly advocated in WW, all the results in WW are
for this special case. 
An important series of papers~\cite{Nemenman:2002fk,nemenman2004entropy,nemenman2008neural} (hereafter abbreviated
as NSB
) considered the generalization where $c = a |Z|$ for any positive constant, $a$. (WW 
is the special case where $a = 1$.)
NSB considered the limit where we have such a $c$, and $|Z|$ is much larger than $N$, the number of samples of $\rho$.
(Therefore, for non-infinitesimal $a$, $c \gg N$.) They showed that in this limit, the samples have little effect on the posterior
moments of $Q$ (e.g., for $Q$ the Shannon entropy). Therefore, the data become irrelevant. In that large $|Z|$ limit, the posterior moments are dominated by the prior over $\rho$
that is specified by $c$. This can be seen as a major shortcoming of WW.

To address this shortcoming, NSB noted that if 
$c$ is a random variable with an associated prior, $P(c)$,
it induces a prior distribution over the values of the Shannon entropy,
$H(\rho) = -\sum_z \rho(z) \ln[\rho(z)]$. Specifically, for the case of a uniform $L$,
for any potential entropy value $h$:
\begin{eqnarray}
P(H = h) &=& \int dcd\rho \; \delta(H(\rho) - h) P(\rho \mid c) P(c) \nonumber \\
&=& \int dcd\rho \; \delta(H(\rho) - h) \frac{\prod_z \rho(z)^{c / |Z| - 1}}{\int d\rho' \; \prod_z \rho'(z)^{c / |Z| - 1}} 
P(c)
\end{eqnarray}
where the integrals over distributions are implicitly restricted to the associated simplex.
NSB then decided that since their goal was to estimate entropy, they would set $P(c)$
so that the prior, $P(H)$, is flat. In essence, they formed a continuous mixture of Dirichlet distributions to (attempt
 to) obtain a flat prior over the ultimate quantity of interest, $H$. The $P(c)$ that results in flat $P(H)$ cannot be written down in closed form, 
but NSB showed how numerical computations can be used to approximate it.

As they are used in practice, both NSB and WW allow $|Z|$ to vary, giving it different
values for different problems, values that are typically set in an \textit{ad hoc} way. (Indeed, if they did not allow 
the number of bins to vary from one
dataset to another, they could not be used on problems with
datasets running over too large a number of bins.) They then set $P(c)$ based on that fixed value of $|Z|$. 

One problem with this is that it 
means the posterior expected value
of the functional, $Q$, can vary depending on how one expresses that functional.
To illustrate this, say that $Z$ is a Cartesian product, $X \times Y$, and let $H(X, Y)$, $H(X)$ and $H(Y)$ refer to the 
entropies of $\rho(x, y) \equiv \rho_{X,Y}(x,y)$, \linebreak $\rho(x) \equiv \rho_X(x) \equiv \sum_y \rho(x, y)$ and 
$\rho(y) \equiv \rho_Y(y) \equiv \sum_x \rho(x, y)$,
respectively. Recall that the mutual information between $X$ and $Y$ under any $\rho$ can be written both as: 
\begin{eqnarray}
I_\rho(X ; Y) &\equiv& -\sum_{x} \rho(x) \ln[\rho(x)] - \sum_{y} \rho(y) \ln[\rho(y)] 
+ \sum_{x, y} \rho(x, y) \ln[\rho(x, y)] \nonumber \\
&=& H(X) + H(Y) - H(X, Y)
\label{eq:first-mut}
\end{eqnarray}
and equivalently as:
\begin{eqnarray}
I_\rho(X ; Y) &\equiv& \sum_{x, y} \rho(x, y) \ln\bigg[\frac{\rho(x, y)}{\rho(x) \rho(y)}\bigg]
\label{eq:second-mut}
\end{eqnarray}

Now, let $\n$ be a set of IID of samples of $\rho_{X,Y}$, and let $\n_X$ and $\n_Y$ be the
associated set of sample values of $\rho_X$ and $\rho_Y$, respectively. 
Then, from Equations~(\ref{eq:first-mut}) and (\ref{eq:second-mut}), the posterior expectation of the mutual
information can be written as either:
\begin{eqnarray}
\E(I(X ; Y) \mid \n) &=& -\E\bigg(\sum_{x} \rho_X(x) \ln[\rho_X(x)] \mid \n_X\bigg) -\E\bigg(\sum_{y} \rho_Y(y) \ln[\rho_Y(y)] \mid \n_Y\bigg) \nonumber \\
&& \;\;\;\;\; \;\;\;\;\; \;\;\;\;\; \;\;\;\;\; \;\;\;\;\; \;\;\;\;\; \;\;\;\;\; \;\;\;\;\; \;\;\;\;\; \;\;\;\;\; + \; \E\bigg(\sum_{x,y} \rho(x,y) \ln[\rho(x,y)] \mid \n\bigg)
 \nonumber \\
&=& \E(H(X) \mid \n_X) + \E(H(Y) \mid \n_Y) - \E(H(X, Y) \mid \n) 
\label{eq:first-post-mut}
\end{eqnarray}
or as:
\begin{eqnarray}
\E(I(X ; Y) \mid \n) &=& \E \bigg( \sum_{x, y} \rho(x, y) \ln\bigg[\frac{\rho(x, y)}{\rho(x) \rho(y)}\bigg] \; \bigg| \; \n \bigg)
\label{eq:second-post-mut}
\end{eqnarray}
(See, for example,~\cite{archer2012Bayesian,archer2013Bayesian}.)

If $P(c)$ is set in a way that depends on the size of the event space, then we would
use a different $P(c)$ to evaluate each of the three terms in Equation~(\ref{eq:first-post-mut}), since the
underlying event spaces ($X$, $Y$ and $X \times Y$, respectively) have different sizes. However, there would
only be a single $P(c)$ used to evaluate the expression in Equation~(\ref{eq:second-post-mut}), the same $P(c)$ as 
used for evaluating the third term in Equation~(\ref{eq:first-post-mut}). \linebreak As a result, under either the NSB or WW approaches,
the values given by Equation~(\ref{eq:first-post-mut}) and Equation (\ref{eq:second-post-mut}) will differ in general;
depending on which definition of mutual information we adopt, we would
get a different estimate of posterior expected mutual information under those approaches.
Indeed, 
to estimate mutual information in the NSB approach, one faces the choice of 
whether to set $P(c)$ to give a uniform prior distribution over values of the mutual information, $P(I(X ; Y))$ (as it would appear, one must, since $I(X ; Y)$ is what one wishes to
estimate), or to set it to give a uniform $P(H(X, Y))$ (as in conventional NSB). It is not clear
how to make this choice, in general.

\section{Contribution of This Paper}

In most of the earlier work on estimating an information-theoretic functional of $\rho$ based
on data, it is assumed that $|Z|$ is fixed, with a few exceptions (see, {e.g.},~\cite{vu2007coverage}).
In many situations, the modeler
is not completely certain \emph{a priori} about the value of $|Z|$, and so should 
treat it as a random variable. For such scenarios, we need to specify a 
joint (hyper)prior, $P(c, |Z|)$, rather than just a prior, $P(c)$. In particular, if we set $c$ from $|Z|$ as in either
WW or NSB, then by specifying our uncertainty in $Z$, $P(|Z|)$, we set the joint prior $P(c, |Z|)$.
Note that for both WW and NSB, this induced $P(c, |Z|)$ is not a product distribution $P(c)P(|Z|)$. (In
particular, in NSB, the distribution, $P(c)$, is set independently of any data, in a way that varies with $|Z|$.)

Jaynes has argued convincingly for setting priors with
invariance arguments concerning the fundamental nature of the problem domain~\cite{jabr03}.
In this paper, we show that the prior, $P(c, |Z|)$, obeys a simple ``Irrelevance of 
Unseen Variables'' (IUV) invariance, 
if and only if $c$ and $|Z|$ are independent, {\em i.e.}, iff $P(c, |Z|) = P(c) P(|Z|)$. 
Therefore, if we require IUV, then rather than specify a full joint distribution over $c$ and $|Z|$, we
only need to specify a distribution over each of $c$ and $|Z|$ separately. This 
greatly reduces the number of degrees of freedom in the prior that we need to specify
(though not as much as would be the case if we used WW or NSB, were we only to
specify $P(|Z|)$).

In this paper, we show that when IUV is obeyed, so that $c$ and $|Z|$ are independent, the value for posterior expected mutual
information does not change depending on whether we use Equation~(\ref{eq:first-mut})
or Equation~(\ref{eq:second-mut}) to define mutual information. 
In proving this, we derive an intermediate result that simplifies the calculation of
some posterior moments. In particular, we show how to use this result to derive the formula for
posterior expected mutual information given in WW in essentially a single line.
We also show that both of these advantages extend to the calculation of multi-information,
one of the ways proposed to generalize mutual information beyond two random
{variables}~\cite{james2011anatomy}. Similarly, since Tsallis entropy with index $q$ is just a weighted
sum over $i$ of the $q$th moments
of the $p_i$, we can evaluate expected Tsallis entropy in closed form using our estimators.  (However, higher-order moments
of the Tsallis entropy do not simplify as easily.)

We then show that when $c$ and $|Z|$ are independent under the prior, and $|Z|$
is averaged over according to a prior $P(|Z|)$ with some reasonable characteristics, the posterior expected
value of \linebreak information-theoretic quantities need not be dominated by the prior. In this sense,
the problem that caused NSB to consider a non-conventional scheme for setting $P(c)$
does not exist if we allow $|Z|$ to be a random variable and require IUV.

We next discuss in detail various fully Bayesian schemes in which the random variables, $c$ and/or $|Z|$,
are integrated over to form estimates of posterior expectations. We also mention some
schemes in which one or the other of those variables is given a single
value (unlike in proper hierarchical Bayes). We run a few computer experiments as cursory ``sanity checks''.
In these, we choose a naive IUV-based hyperprior and compare the associated estimators of posterior expected entropy
and of posterior mutual information
to the estimators considered {in}~\cite{Nemenman:2002fk, nemenman2011coincidences, vu2007coverage}.
We find that the IUV-based estimator performs quite favorably.


There are several subtleties in how one models the statistical generation of $Z$, issues
that do not arise if $Z$ is fixed ahead of time. One
of them involves the mapping of each newly sampled draw of $\rho$ to an element of $\n$,
{\em i.e.}, to a label for that draw. To formally justify the ``intuitively obvious'' model of how $Z$
is generated that we have used up to now, we describe in detail a mapping of draws of $\rho$
to elements of $\n$ that justifies that model.

After this, we describe a change one might make to the model of
how $Z$ is generated that would appear to be innocuous. We show that, in fact, this change can substantially
affect the resultant estimations. Concretely, say there is a space, $\hat{Z}$, that is 
a grid of photoreceptors and that $\rho$ is a distribution over $\hat{Z}$ that is IID sampled 
to generate counts of photons that are reflected from an object and focused onto elements of $\hat{Z}$. Say we know that
the object being imaged may be occluded, so that, in fact, only a subset, $Z \subseteq \hat{Z}$,
of the grid points can have a nonzero probability of a photon count. However, say we are uncertain 
of the size of $Z$ and, therefore, of which precise pixels in $\hat{Z}$ it corresponds to. We show that the value
of the posterior expected entropy in this scenario is different from its value in the conventional
scenario, in which we are also uncertain of $Z$'s size, but there is no encompassing
set $\hat{Z}$ from which $Z$ is formed.

This touches on the more general issue of the epistemological foundation of the probabilities (and
probabilities of probabilities) considered in this paper. This issue, involving concepts like ``degree of belief'' and ``objective
probability'', is deep and quite important, being fundamental to
the differences between Bayesian and sampling theory
statistics (see~\cite{wolp96c, wolp95b} for a discussion). 

Here, we do not
grapple with this issue. Rather, we adopt the ``pragmatic Bayesian'' perspective implicit in all earlier Bayesian
work on the problem of estimating information-theoretic quantities from samples (including
WW and NSB, in particular) and simply use probabilities as a part of a self-consistent calculus of uncertainty. Bayesian reasoning often relies on a choice of prior, and the work presented here is no exception. We emphasize that that is no ``one true prior''; rather, the statistician must match the choice of model to their own prior knowledge about the system. While we have done our best to chose a set of priors with generality sufficient to avoid some of the common biases identified in the past, one of the main goals of this paper is to present our results, so that the readers may adapt our methods to the particular nature of their own research.

Some of the experiments reported here were run using the publicly available package, Thoth, available at {http://thoth-python.org}
. The reader is directed to Appendix B for associated proofs. Note that as discussed
in WW, much of the analysis below can be modified for inference of arbitrary functionals
of $\rho$ from IID samples of $\rho$ (e.g., estimation of covariances from IID samples
rather than mutual information). The analysis is not limited to inferring information-theoretic
functionals from IID samples.
%

\section{Preliminaries}
\label{sec:prelim}

For any finite space, $U$,
we write $\Delta_U$ to mean the simplex of possible distributions over $U$.
We will also write $T_k$ to mean the set of all $k$-dimensional vectors whose
components are all integers greater than~zero.

Throughout this paper, we restricted attention to Dirichlet distributions
over distributions $\rho \in \Delta_Z$. For~a fixed $c$ and $Z$, we write such a distribution as:
\begin{eqnarray}
\D_{c,Z}(\rho) &=& \frac{ \prod_z \rho(z)^{[c / |Z|] - 1}} {\int d\rho' \; \prod_z \rho'(z)^{[c / |Z|] - 1}}
\label{eq:dirichlet}
\end{eqnarray}
where we require $c$ to be non-negative.

Say we are given a dataset, $\n$, of counts for each of the elements of $Z$,
where $N \equiv \sum_z n_z$. 
For the Dirichlet prior, the posterior distribution is:
\begin{eqnarray}
P(\rho \mid \n, c, |Z|) &=& \frac{\D_{c,Z}(\rho) \prod_z \rho_z^{n_z}} {\int d\rho' \; \D_{c,Z}(\rho') \prod_z (\rho'_z)^{n_z}} \nonumber \\
&\equiv& \frac{\prod_z \rho_z^{n_z -1 + c / |Z|}} {G(\n, c, |Z|)}
\label{eq:post-rho}
\end{eqnarray}
Using~\cite{wowo95a}, we can calculate:
\begin{eqnarray}
G(\n, c, |Z|) &=& \frac{\prod_z \Gamma(n_z + c / |Z|)} {\Gamma(N + c)}
\end{eqnarray}
We will sometimes write the posterior given by Equation~(\ref{eq:post-rho}) as $\D_{c,Z}(\rho \mid \n)$.

Note that $Z$ is implicitly a random variable in Equation~(\ref{eq:post-rho}), since we condition on $|Z|$ there. This means that
the event space over which $\rho$ is defined is also a random variable, as is the event space
over which $\vec{n}$ is defined. This means that when we average over $Z$'s below, we 
must take a bit of care to define the space of all $\rho$'s and the space of all $\vec{n}$'s, since $\rho$ can be an element of
any finite unit simplex and similarly for $\vec{n}$. 

When this issue arises, we will
define the set of all triples of $Z$, $\n$ and $\rho$
as the infinite union of the triples, $(Z_1, \Delta_{Z_1}, T_{|Z_1|})$, $(Z_2, \Delta_{Z_2}, T_{|Z_2|})$, $(Z_3, \Delta_{Z_3}, T_{|Z_3|})$, \textit{etc}.,
where each $Z_i$ is defined as the set, $\{1, \ldots, i\}$. In such a fully formal approach, joint probability distributions over 
$Z, \n$ and $\rho$ are defined over that infinite union by:
\begin{eqnarray}
P(Z, \rho, \n) &=& 0 {\mbox{ unless }} \rho \in \Delta_Z, \n \in T_{|Z|} 
\label{eq:var_Z_formal}
\end{eqnarray}
and then using the multinomial and Dirichlet distributions to define the values of the conditional
distributions, $P(\n \mid \rho, Z)$ and $P(\rho \mid Z)$, respectively, when the condition
in Equation~(\ref{eq:var_Z_formal}) is obeyed.
For the simplicity of the exposition, we will minimize our use of this fully formal approach here.


We define $\SSS(\n)$ to be the support of the dataset, $\n$, within $Z$. We write $I(.)$ to be 1/0, depending
on whether its logical expression argument is true/false. 
For the case where $Z = X \times Y$, we define $\n_X(x) \equiv \sum_y \n(x, y)$
and $\rho_X \equiv \sum_y \rho(x, y)$. 
For use below, as in WW, we define $\Delta\Phi^{(1)}(z_1,z_2) \equiv \Psi^{(0)}(z_1)-\Psi^{(0)}(z_2)$, where 
$\Psi^{(0)}$ is defined as $d[\ln{\Gamma{(z)}]} / dz$.

\section{Irrelevance of Unseen Variables}
\label{sec:iuv}
\vspace{-12pt}

\subsection{The Problem of Unseen Variables}

In general,
there may be an ``unrecorded'' or ``hidden'' variable, $y \in Y$,
whose values are not recorded in our dataset of values, $x \in X$.
As an example, say our data are the set of all changes
in the value of the US stock-market between opening and
closing on all days it was open since 1970.
A hidden variable is the age of the person recording each of the measurements
when they made the recording. The change in stock market
value is $X$, and the age of the person recording the value is $Y$.

The hidden variable in this example is chosen to be extreme, in that
knowledge of its existence is clearly irrelevant to the statistical estimation
problem. However it is hard to imagine an estimation problem where there
are \emph{not} in fact a potentially infinite number of such hidden variables. Some 
could perhaps be dismissed, as ``clearly irrelevant, and therefore,
not worthy of consideration''. However, this approach is hard to justify axiomatically,
since there are, of course, many instances when hidden variables \emph{may} have
some relevance to the estimation problem.

This issue is a bit of a philosophical hornet's nest. If at all possible, we would
like to avoid having to consider it. In fact, we can do this by requiring that
the existence of any hidden variables has no effect on our Bayesian estimate.
More precisely, one can require that whether one does or does not have a single unseen random variable in
one's model, and the (finite) size of its event space, if it does exist, must have no impact on posterior expected values of (functionals
of the distribution over) the seen variables~\cite{davidone}. How to do this is the subject of this section. In the following section, we illustrate
the practical benefits that result.

To analyze scenarios involving hidden variables, write $Z = X \times Y$, and
say we have recorded the dataset, 
$\n_X(x) \equiv \sum_y \n_{X,Y}(x, y)$, not the full dataset, $\n_{X,Y}$.
Next, write $\rho_Z$ as a matrix
of real numbers, $\{\rho_{X, Y}(x, y) : x \in X, y \in Y\}$. We can re-express any such $\rho_{X, Y} \in \Delta_{X \times Y}$ 
in an alternative coordinate system, as $(\rho_X, \rho_{Y \mid X})$,
where $\rho_X(x) \equiv \sum_y \rho_{X, Y}(x, y)$, and
$\rho_{Y \mid X}$ is the set of $|X||Y|$ real numbers given by $\rho(x, y) / \rho_X(x)$ for all
$x \in X, y \in Y$. (Note that under a Dirichlet prior, no matter what our data are,
there is both zero prior probability and zero posterior probability of a $\rho_X$, such that $P(\rho_X)(x) = 0$
for some $x$.) Therefore, for all $x, y$, $\rho_{X,Y}(x, y) = \rho_{Y \mid X}(y \mid x) \rho_X(x)$.

If we allow for
such hidden variables, $Y$, then to do a proper Bayesian analysis, we must
specify a prior, $P(\rho_{X,Y})$ (or just $P(\rho)$, for short) on the space, $\Delta_{X \times Y}$, {\em i.e.}, on the space
of $\rho$'s that run over $X \times Y$. It does not suffice to specify
just a prior, $P(\rho_X)$, defined over $\Delta_{X}$. Moreover, axiomatic
derivations of Bayesian analysis counsel us to set
this prior \emph{without concern for what our likelihood function
will be}. (The prior is our model of the underlying physical system.
The likelihood instead has to do
with the observation apparatus we happen to have handy to observe that system.)
This implies that $P(\rho)$ should not reflect the fact that $X$ is observed
and $Y$ is not, since what variable is observed is determined by the likelihood.
Therefore, in particular, if we set $P(\rho)$ based on the size of the underlying event space,
it should be based solely on the size of the space, $X \times Y$,
with no consideration for just the size of $X$.
%


Now, in general, we do not even
know how many values of a hidden variable $Y$ there are, simply
that there may (!
) be some. Due to this uncertainty, we must let the
cardinality of the set of hidden variables ``float'' as a random variable.
Moreover, very often we are interested in a functional $Q(\rho^X(x)) = Q(\sum_y \rho_{X, Y}(x, y))$
of the observed variable, and this functional can be anything, depending
on the statistical question we are interested in. 

One might worry that for some such functional, $Q$, a Dirichlet
prior over the the joint observed and hidden variables, $P(\rho_{X,Y})$,
and some (visible) data vector, $\n_X$, the associated value of
$\E(Q \mid \n_X)$ would vary depending on the number
of degrees of freedom of the hidden variable, $|Y|$. If this
were the case, how we set the prior
over the number of degrees of freedom of the hidden variable would
matter, which would confront us with the problem of how to set it.
This would seem to be an intractable problem, since, in general,
there are an infinite number of choices for what the hidden variable is.

The desideratum analyzed in this paper is that this problem does not
arise. Formally, this is equivalent to requiring that no matter what $Y$ is,
the associated value of $\E(Q \mid \n_X)$ is exactly what it would be if
there were no hidden variable at all:
\begin{definition} 
A distribution $\pi(c, |Z|)$ obeys \textbf{Irrelevance of Unseen Variables} (IUV) iff,
for all finite spaces, $X$ and $Y$, data vectors, $\n_X$, and functions, $Q$, defined over $\rho_X$:
\begin{eqnarray}
\int dc \; \pi(c \; | \; |X|) \int d\rho_X \; Q(\rho_X) \D_{c, X}(\rho_X \; | \; \n_X) && \nonumber \\
&& \!\!\!\!\!\!\!\!\!\!\!\!\!\!\!\!\!\!\!\!\!\!\!\!\!\!\!\!\!\!\!\!\!\!\!\!\!\!\!\!= \;
 \int dc \; \pi(c \mid |X| \times |Y|) \int d\rho_{X,Y}\; Q(\rho_X) \D_{c, X \times Y}(\rho_{X,Y} \; | \; \n_X) \nonumber 
\label{def:iuv-def}
\end{eqnarray}
\end{definition}

To help understand this desideratum, note that uncertainty about the size of a hidden variable space, $|Y|$, is different
from uncertainty about the size of the observed variable space, $|X|$.
Indeed, one could argue that $|Y|$ is always essentially infinite,
up to any kinds of limits imposed by quantum mechanics. \linebreak 
(As an illustration for the stock-market example, in addition to the age of the
recorder of the \linebreak stock-market's change in value, all other
characteristics of the recorder are ``hidden'' and, therefore, arguably
should be included in $Y$.) Moreover, in general,
as the number of (IID) data grows, we will get more certain about $|X|$,
or at least about the number of $x$ for which $\rho(x)$ exceeds some preset
threshold. In contrast, the size of the dataset has no effect on our
uncertainty concerning $|Y|$; the latter is purely prior-dominated. 
Both of these properties mean that statistically
estimating $|Y|$ is a more fraught exercise
than estimating $|X|$.

Our desideratum says that $\pi$ is arranged so that these difficulties 
are irrelevant. For a $\pi$ obeying IUV, uncertainty in the value of $|Y|$ has no effect on our estimate of
a functional $Q(\rho_X)$.
In contrast, uncertainty in $|X|$ has a major effect on all estimators of functionals
$Q(\rho_X)$ that we know of (including the ones we introduce below), regardless of $\pi$.
Indeed, how to estimate the size of $X$ from the observed data is so important
that it has been analyzed for decades in statistics, under the name ``coverage 
estimators''~\cite{bunge1993estimating,vu2007coverage}. 

The $\pi$'s used in both NSB and WW violate IUV. This is at the heart of their problems
in estimating mutual information, which were discussed in Section~\ref{sec:background}.

\subsection{Dirichlet-Independent (DI) Hyperpriors}

Let $T$ be a partition of $Z$. Then, there is a map, $K$, taking any distribution, $\rho$, over $Z$
to a distribution, $\rho_T$, over the elements, $t \in T$. 
Using $K$, any distribution over $\rho$'s induces a distribution over $\rho_T$'s.
In particular, Dirichlet distributions have the very nice property that the distribution $\D_{c,Z}(\rho)$
over $\rho$'s induces the distribution, $\D_{c,T}(\rho_T)$, over $\rho_T$'s, where the baseline
distribution for $T$ is given by applying $K$ to the baseline distribution over $Z$. The crucial point about this property for us
is that there is the same concentration parameter, $c$, in both Dirichlet distributions;
Dirichlet distributions are consistent under marginalization. (Indeed, this
property serves as a common definition of Dirichlet processes, the extension of
Dirichlet priors to infinite spaces.)

As a special case, if $Z = X \times Y$, then $X$ specifies a partition of
$Z$ in which each partition element is of the form $\{(x, y) : y \in Y\}$ for a different $x \in X$. Therefore, a 
Dirichlet distribution generating $\rho$'s over $Z$ induces
a Dirichlet distribution generating $\rho_X$'s over $X$ that has the same concentration parameter.

Now note that because we are using a
Dirichlet prior, the posterior, $P(\rho \mid \n_X)$, is a Dirichlet distribution. 
In light of the marginalization consistency of Dirichlet distributions just described, this means that the induced
posterior, $P(\rho_X \mid \n_X)$, will also be a Dirichlet distribution,
\emph{with the same value of $c$}. This suggests that the posterior
expectation of any functional of $\rho_X$ will be the same,
whether we evaluate it in $X$ or $X \times Y$, so long as $c$
is the same for both evaluations.

We can formalize this with the following lemma, proven in Appendix B:

\begin{lemma}
Fix any finite spaces, $X$ and $Y$, any set, $\n_{X,Y}$, of counts of each of the elements in $X \times Y$ and 
any two concentration parameters, $c$ and $c'$. Then, $c = c'$ iff:
\begin{eqnarray}
 \int d\rho_X \; Q(\rho_X) \D_{c,X}(\rho_X \; | \; \n_{X}) &= &
 \int d\rho_{X,Y} \; Q(\rho_X) \D_{c',X\times Y} (\rho_{X,Y} \;| \; \n_{X,Y}) \nonumber
\end{eqnarray}
for all functions, $Q$, defined over $\Delta_X$.
\label{lemma:1}
\end{lemma}

There are several noteworthy implications of Lemma~\ref{lemma:1}. To see the first one,
consider the following modification of the definition of IUV, which involves $\n_{X,Y}$,
the count vector of both seen and unseen~bins:
\begin{definition}
A distribution, $\pi(c, |Z|)$, obeys \textbf{strengthened IUV} iff
for all finite spaces, $X$ and $Y$, data vectors, $\n_{X,Y}$, and functions, $Q$, defined over $\rho_X$:
\begin{eqnarray}
\int dc \; \pi(c \; | \; |X|) \int d\rho_X \; Q(\rho_X) \D_{c, X}(\rho_X \; | \; \n_X) &= & \nonumber \\
&& \!\!\!\!\!\!\!\!\!\!\!\!\!\!\!\!\!\!\!\!\!\!\!\!\!\!\!\!\!\!\!\!\!\!\!\!\!\!\!\!
 \int dc \; \pi(c \mid |X| \times |Y|) \int d\rho_{X,Y}\; Q(\rho_X) \D_{c, X \times Y}(\rho_{X,Y} \; | \; \n_{X,Y}). \nonumber 
\end{eqnarray}
\label{def:iuv-strong}
\end{definition}
An immediate corollary of Lemma~\ref{lemma:1} is the following:
\begin{corollary}
IUV implies strengthened IUV.
\label{coroll:iuv-strong}
\end{corollary}
{\noindent \textbf{Proof:}}
The integral on the RHS
 in the equation in Lemma~\ref{lemma:1} is the inner integral on the RHS
of Definition~\ref{def:iuv-strong}. In~addition, the integral on the LHS
 in the equation in Lemma~\ref{lemma:1}
is the inner integral on the RHS of Definition~\ref{def:iuv-def}. Therefore, applying 
Lemma~\ref{lemma:1} for $c = c'$ establishes the corollary.
\qed

$ $

%
%
%
\noindent Corollary~\ref{coroll:iuv-strong} means that if IUV holds, then it does not matter how the counts, $\n_{X,Y}(x, y)$, are
apportioned over $Y$, as far as calculating the associated posterior
expected value of $Q$ is concerned.

Assume there are no unobserved components of our data.
Then, by Corollary~\ref{coroll:iuv-strong}, for any functional, $Q$, that depends purely on $\rho_X$,
if IUV holds, we can evaluate expected
$Q(\rho_X)$ conditioned on $\n_{X,Y}$ (given by an integral over $\Delta_{X \times Y}$)
by calculating expected $Q(\rho_X)$ conditioned on $\n_X$ (given by an integral over $\Delta_X$).
In Section~\ref{sec:bennies} below, we show that this property of IUV substantially simplifies calculations of expected
moments of functionals defined over multi-dimensional spaces,
e.g., the calculation of posterior expected mutual information.

To establish a second implication of Lemma~\ref{lemma:1},
multiply both sides of the equality it establishes by $\sum _{\n_{X,Y}} P(\n_{X,Y} \mid \n_X)$. That
leaves the LHS of that equality unchanged. However, it changes the RHS to
$ \int d\rho_{X,Y} \; Q(\rho_X) \D_{c',X\times Y} (\rho_{X,Y} \;| \; \n_{X})$. This
provides the following corollary:
\begin{corollary}
Fix any finite spaces, $X$ and $Y$, 
concentration parameter $c$, function $Q$ defined over $\Delta_X$ 
and set $\n_{X,Y}$ of counts of each of the elements in $X \times Y$. 
Then:
\begin{eqnarray}
 \int d\rho_X \; Q(\rho_X) \D_{c,X}(\rho_X \; | \; \n_{X}) &= &
 \int d\rho_{X,Y} \; Q(\rho_X) \D_{c,X\times Y} (\rho_{X,Y} \;| \; \n_{X}) \nonumber
\end{eqnarray}
\label{coroll:partial-iuv}
\end{corollary}
\noindent The integral on the RHS of Corollary~\ref{coroll:partial-iuv} is 
the inner integral on the RHS of the definition of IUV, Definition~\ref{def:iuv-def}. 
This establishes that
if the conditional prior, $\pi(c \mid |X|)$, equals the conditional prior, $\pi(c \mid |X||Y|)$, then IUV holds.
%

We will use the term \textbf{Dirichlet-independent hyperprior} (DI hyperprior)
to refer to any prior of the form $\pi(c, |Z|) = \pi(c) \pi(|Z|)$
over the hyperparameters of a Dirichlet distribution with uniform baseline distribution.
We now present the main result of this section, which is that if we restrict
attention to hyper priors meeting a particular technical 
condition~\cite{davidtwo}, then IUV and a DI hyperprior are equivalent, as proven in Appendix B:
\begin{proposition}
\label{prop:iuv-weak}
Assume that for any two finite spaces, $X and Y$, and associated count vector $\n_X$, there
exists an $\epsilon > 0$, such that:
\begin{eqnarray*}
((1 + \epsilon) |X|)^c \frac{\pi(c \mid |X|) - \pi(c \mid |X||Y|)} {C(c, \n_X) }
\end{eqnarray*}
is infinitely differentiable with respect to $c$ at $c = 1$ and that its Fourier transform is analytic. Then,
$\pi$ is a DI hyperprior $\Leftrightarrow$ IUV holds.
\end{proposition}

\noindent We can combine these results to see that under the conditions in Proposition~\ref{prop:iuv-weak},
we have a DI hyperprior iff strengthened IUV holds.

In some situations, the number of possible values of the ``hidden variable''
will vary with $x \in X$. This is a generalization of the currently considered
scenario, where rather than $Z = \cup_{x \in X} Y_x$, where $|Y_x|$ is the
same for all $x$, we allow the sizes, $|Y_x|$, to vary with $x$. 
One can modify the definitions of IUV and strengthened IUV given above to address
this situation, and all the results presented above (appropriately modified) still hold.
This is due to the fact that Dirichlet distributions are ``consistent under
marginalization'', as described at the beginning of this section.

\section{Calculational Benefits}
\label{sec:bennies}

In addition to its ``theoretical'' advantage of being equivalent to a natural desideratum (IUV), DI hyperpriors also have
practical advantages.
In particular, recall that WW derived a complicated expression for posterior expected mutual information based on 
Equation~(\ref{eq:second-post-mut}). Part of what made that expression complicated
was that it used a posterior over $\Delta_{X,Y}$ to evaluate $H(X)$ and $H(Y)$,
as well as $H(X, Y)$.

However, when we have a DI hyperprior, the posterior expectation of $H(X)$ 
conditioned on $\n_X$ is the same, whether
we evaluate it using a posterior over $\Delta_X$ or
evaluate it using a posterior over $\Delta_{X,Y}$. This means we can evaluate that posterior expected
mutual information by summing the posterior expected $H(X)$ under a posterior over $\Delta_X$, 
posterior expected $H(Y)$ under a posterior over $\rho_Y$ and posterior expected $H(X, Y)$ under 
a posterior over $\Delta_{X,Y}$. In turn, each of those three expectation values are given by a relatively
simple formula from WW [Equation~(\ref{eq:post_ent_mean}) in Appendix A].

This property concerns a scenario where we only have data $\n_X$. However, often when we want 
to estimate mutual information, we will have the full dataset with no unseen components, $\n_{X,Y}$.
For that case, we can use Lemma~\ref{lemma:1} to justify decomposing the
posterior expected mutual information into a sum of three posterior expected entropies and then
evaluate those posterior expected entropies separately using 
Equation~(\ref{eq:post_ent_mean}) in Appendix A. This directly gives us the following result, which
does not rely on~IUV:
\begin{corollary}
Let $X$ and $Y$ be two spaces and $\vec{n}$, a sample generated by IID generating a
distribution, $\rho$, across $X \times Y$, where $\rho$ was generated by sampling a Dirichlet
prior with concentration parameter $c$. Then:
\begin{eqnarray*}
\E(I(X ; Y) \mid {\vec{n}}, c) &=& \sum_{x,y} \Delta \Phi^{(1)} (\n(x,y) + 1 + c / |X||Y|, N + |X||Y| + 1) 
			\frac{\n(x,y) + c /|X||Y|}{N + c} \nonumber \\
 && - \sum_{x} \Delta \Phi^{(1)} (n_X(x) + 1 + c / |X|, N + |X| + 1) 
			\frac{ n_X(x) + c /|X|}{N + c} \nonumber \\
 && - \sum_{y} \Delta \Phi^{(1)} (n_Y(y) + 1 + c / |Y|, N + |Y| + 1) 
			\frac{ n_Y(y) + c /|Y|}{N + c} \nonumber \\
\end{eqnarray*}
\end{corollary}
\noindent By Corollary~\ref{coroll:iuv-strong}, the analogous equality holds when we marginalize over $c$
rather than condition on it, so long as we assume IUV. 

Similarly to how the posterior first moment of mutual information can
be defined using either Equation~(\ref{eq:first-mut})
or Equation~(\ref{eq:second-mut}), so can the posterior second moment:
\begin{eqnarray}
E(I^2 \mid \n) &=& \E(H(X)^2 \mid \n_X) + \E(H^2_Y \mid \n_Y) + \E(H^2_{X,Y} \mid \n) \nonumber \\
&& \;\;\;\;\;\;\;\;\;\;\;- \; 2 \E(H(X) H(Y) \mid \n) - 2\E(H(X) H(X, Y) \mid \n)
-2 \E(H(Y) H(X, Y) \mid \n) \nonumber \\
\label{eq:first-second-post-moment-mut}
\end{eqnarray}
or, alternatively:
\begin{eqnarray}
\E(I^2 \mid \n) &=& \E \bigg( \bigg( \sum_{x, y} \rho(x, y) \ln\bigg[\frac{\rho(x, y)}{\rho(x) \rho(y)}\bigg] \bigg)^2 \; \bigg| \; \n \bigg)
\label{eq:second-second-post-moment-mut}
\end{eqnarray}
where the RHS of Equation~(\ref{eq:second-second-post-moment-mut}) is evaluated using a Dirichlet prior over $\Delta_{X,Y}$, while
the first two terms on the RHS in Equation~(\ref{eq:first-second-post-moment-mut}) are instead evaluated
using a Dirichlet prior over $\Delta_{X}$ and over $\Delta_Y$, respectively.

Using the DI hyperprior, one gets the same answer whichever one of these expansions one uses. In~addition,
the first three terms in Equation~(\ref{eq:first-second-post-moment-mut}) can be simplified under
the DI hyperprior and evaluated using the formula in WW for posterior expected entropy. 
Unfortunately though, the remaining terms, e.g., $\E(H(X) H(X, Y)) \mid \n)$, cannot
be simplified the same way; evaluating them seems to require the kinds of techniques used in WW to evaluate the posterior variance
of mutual information. This also applies to the use of our estimators for the computation of Tsallis entropy, since they also require use of higher-order moments..

There have been many ways proposed to generalize the idea of mutual information
beyond two random variables. One of the most prominent is the \textbf{multi-information}
of a set of random variable (see, {e.g.}, Ref.~\cite{james2011anatomy}).
Just like mutual information among a pair of random variables can be defined either as a sum of the entropies of
subsets of those random variables or as a function over the full event space, the same is
true of multi-information. The definition of multi-information in terms of the entropies of subsets of the random variables
is:
\begin{eqnarray}
I(X_1, X_2, \ldots) &\equiv& \sum_i H(X_i) - H(X_1, X_2, \ldots)
\end{eqnarray}

Just as requiring IUV means that we do not have to worry about which of the ways
to express the mutual information of two random variables (when calculating the posterior expectated
mutual information based on data concerning only one of those random variables), the same is true of multi-information;
requiring IUV means that we do not have
to worry about how to express multi-information among a set of random variables
when calculating its posterior expectation based on
data concerning only a subset of those random variables.
Moreover, just as IUV greatly simplifies the calculation of posterior expected
information when our data concerns all the random variables, allowing us to
just repeatedly apply Equation~(\ref{eq:post_ent_mean}) in Appendix A, the same is true with multi-information.

\section{Uncertainty in the Concentration Parameter and Event Space Size}
\label{sec:uncertainties}

Even when one adopts a DI hyperprior, so that prior $c$ and prior $|Z|$ are
statistically independent, there is still the issue
of how to set $P(c)$ and $P(|Z|)$. In this section, we discuss some aspects of this
issue.

\subsection{Uncertainty in $c$}
\label{sec:single_c}

There are several natural choices of $P(c)$. For example, if we view $c$ as a scale
parameter, a~logarithmic prior ($P(c) \propto 1/c$ up to a very large cut-off in $c$) would be reasonable. 

%
%

Another approach is to set $c$ to a single value that is ``optimal'' in some sense.
Grassberger argued for setting $c = 1$ based on minimizing a rough approximation to the
statistical bias. In addition, as~subsequently pointed out by NSB, for a fixed $Z$,
the choice of $c$ equal to unity gives near-maximal prior variance of the entropy, {\em i.e.}:
\begin{equation}
\E( [H - \E(H \mid c, Z)]^2 \mid c, Z)
\; =\; \frac{c/|Z|+1}{c+1}\psi_1(c/|Z|+1)-\psi_1(c+1)
\end{equation}
is near its maximum at $c = 1$. Confirming this nice property of $c = 1$,
numerical calculation finds that, as $|Z|$ goes to infinity, the value of $c$
maximizing that prior variance is $c_{\mathrm{max}}\approx0.9222$. For small numbers of bins, $c_\mathrm{max}$ is smaller ({e.g.}, for $|Z|$ equal to 5, $c_{\mathrm{max}}\approx0.6997$).

Of course, one could also set $c$ via a scheme other than hierarchical Bayes. In particular, setting it via maximum
likelihood ({\em i.e.}, ML-II~\cite{berg85}) should often be reasonable.

\subsection{Likelihood of the Event Space Size}
\label{sec:|Z|_like}

Recalling the definition of $G$ from Section~\ref{sec:prelim}, we can write:
\begin{eqnarray}
P(\n \mid c, |Z|) &=& \frac{G(\n, c, |Z|)} {\sum_\n G(\n, c, |Z|)} 
\end{eqnarray}
Using the results of WW to evaluate this, we get:
\begin{eqnarray}
P(\vec{n} \mid c, |Z|) &=& \frac{\Gamma(c)} {\Gamma(c/|Z|)^{|\SSS(\n)|}} \frac{\prod_{i=1}^{|\SSS(\n)|} \Gamma{(n_i+c/|Z|)}} {\Gamma{(N+c)}}
\label{post1}
\end{eqnarray}
Note that $\Gamma(x)$ diverges as $1/x$ as $x\rightarrow0$. Thus, the factor of $\Gamma(c/|Z|)^{|\SSS(\n)|}$ in the denominator means that 
Equation~(\ref{post1}) is strongly weighted towards $|Z|$ small (but strictly greater than $|\SSS(\n)|$, the number of observed bins.)

\begin{figure}[H]\centering
\includegraphics[width=4.2in]{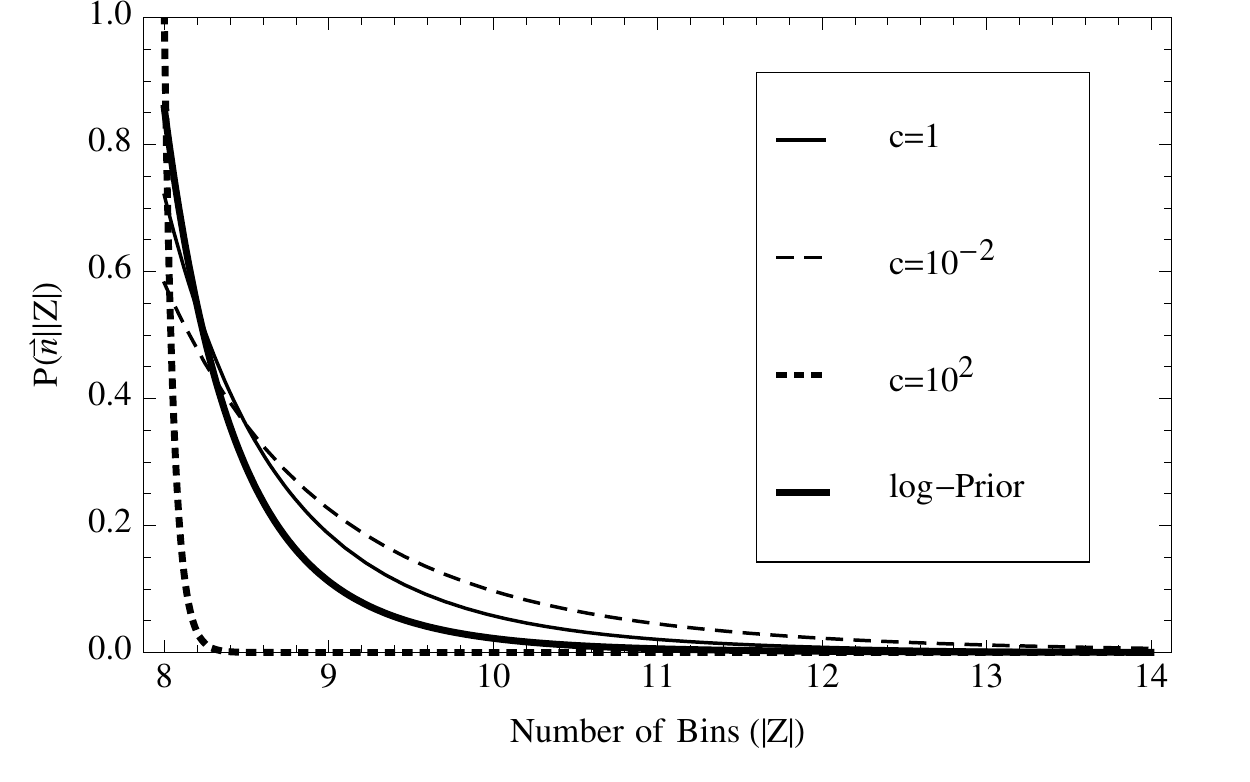}
\caption{Likelihoods for a dataset of one thousand samples drawn from a distribution $\rho$ that
was, in turn, drawn from a Dirichlet prior. The concentration parameter was $c = 1$, and $|Z|$ was 100. (The actual dataset was $\{691, 232, 24, 17, 14, 10, 6, 6\}$,
with the remaining 92 entries equaling zero.)
{\bf Thick solid line}: $P(\vec{n} \mid |Z|)$ with a logarithmic prior for $c$ [Equation~(\ref{int})]. Mean value: $8.2$. {\bf Thin solid line}: $P(\vec{n} \mid
|Z|,c)$ with $c = 1$ (mean value: $8.4$). {\bf Dashed line}: $c = 0.01$ (mean value: $8.8$). {\bf Dotted line}: $c = 100$ (mean value $8.0$.) Note that the maximum likelihood (ML) value is always $|\SSS(\n)|$, the number of observed bins.}
\label{est}\vspace{-12pt}
\end{figure}

We have fixed the constant, $c$, in Equation~(\ref{post1}). We could integrate over it instead, getting:
\begin{equation}
P(\vec{n} \mid {|Z|})=\int \frac{\Gamma(c)}{\Gamma(c/{|Z|})^M}\frac{\prod_{i=1}^M \Gamma{(n_i+c/{|Z|})}}{\Gamma{(N+c)}}P(c)~dc
\label{int}
\end{equation}
where $P(c)$ must be independent of ${|Z|}$ in order to preserve IUV.
Choosing different values of $c$ shifts the posterior distribution, as can be seen in Figure~\ref{est}.

Figure~\ref{est} plots the likelihoods, $P(\n \mid |Z|)$ and $P(\n \mid c, |Z|)$, for various
values of $c$. Note that $\SSS(\n)$ is far smaller than $|Z|$.
This reflects the fact that for $c \ll |Z|$, $\rho$
is likely to be highly peaked about only a few bins and close to zero for
all others.
Note how strongly this likelihood prefers small values of $|Z|$. 
Even when the likelihood is evaluated with the value of $c$ that was used to 
generate $\n$, it still strongly prefers $|Z|$ values that are far smaller
than the one that was used to generate $\n$~\cite{davidthree}.

\begin{figure}[H]\centering
\includegraphics[width=5in]{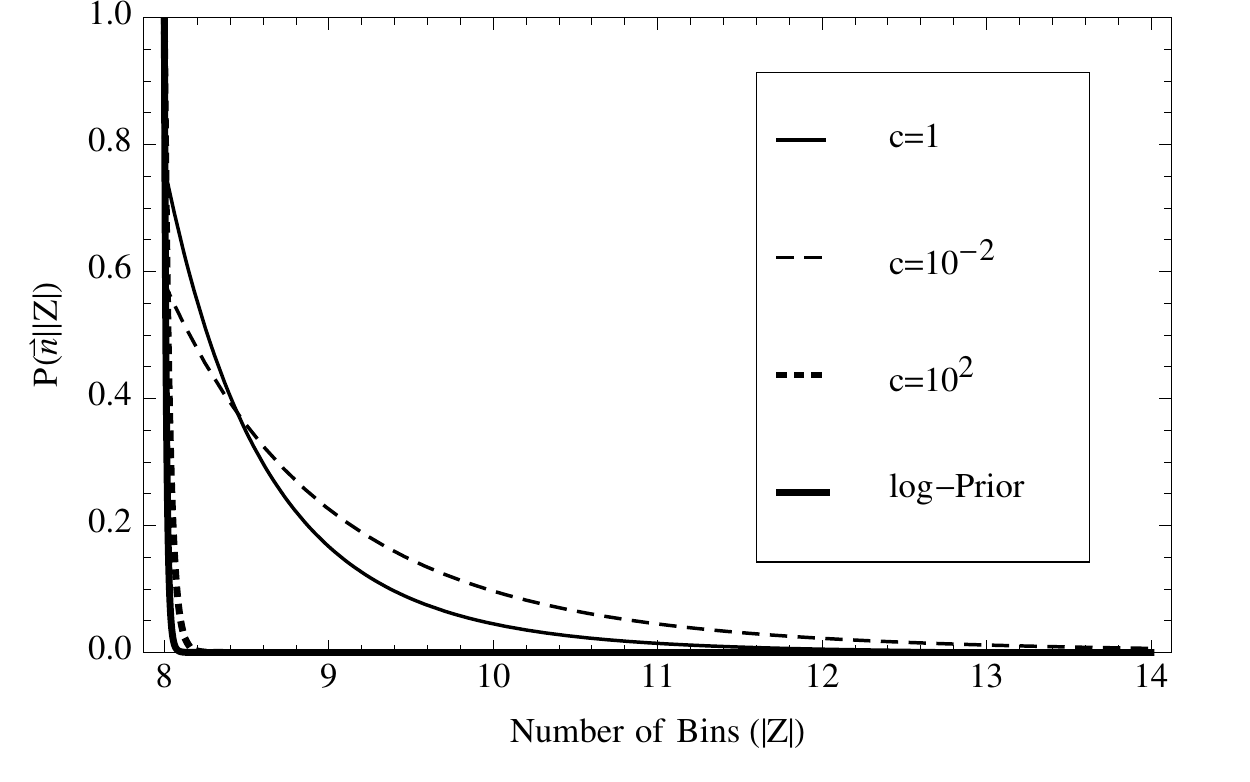}
\caption{Likelihoods for a dataset consisting of eight bins with 125 counts each, with the remaining 92 entries being zero.
{\bf Thick solid line}: $P(\vec{n} \mid |Z|)$ with logarithmic prior for $c$ (Equation~\eqref{int}.) Mean value: $8.0$. {\bf Thin solid line}: $P(\vec{n} \mid
|Z|,c)$ with $c = 1$ (mean value: $8.3$). {\bf Dashed line}: $c = 0.01$ (mean value: $8.8$). {\bf Dotted line}: $c = 100$ (mean value $8.0$). Note that the ML value is always $|\SSS(\n)|$, the number of observed bins.}
\label{unif-est}
\end{figure}

In {Figure}~\ref{unif-est}, we again evaluate $P(\n \mid |Z|)$ and $P(\n \mid c, |Z|)$,
where $\SSS(\n)$ again equals eight, but now, $\n$ is uniform with the value 125 over the eight occupied
bins. This dataset again implies with a high probability that $\rho$ was highly peaked
about the eight occupied bins and close to zero elsewhere. Therefore, again, the likelihood
has a strong preference towards small values of $|Z|$. 

To understand these results, recall ``Bayes factors'', which arise in Bayesian inference of the dimension
of an underlying stochastic model based on samples of that model.
These factors cause a strong \emph{a priori} preference of 
the model dimension's likelihood towards small values.
The preference of $P(\n \mid |Z|)$ for small $|Z|$
is a similar phenomenon, with the size of the space $Z$ playing an analogous
role here to the model dimension in Bayes factors. In both cases,
it is the greatly increased likelihood of the data for a distribution from
the smaller model that causes that model to have greater likelihood.

By comparing Figures~\ref{est} {and}~\ref{unif-est}, we see that for a logarithmic prior (one that
is agnostic about $c$), changing the data without changing $|Z|$ or $\SSS(\n)$ can
have a marked effect on the likelihood, even when $|Z| \gg |\SSS(\n)|$~\cite{davidfour}.
Evidently, then, so long as we integrate over $c$'s for each $|Z|$ (as in NSB) rather than fix a single value (as
in WW), the dependence of the likelihood of $|Z|$ on the data is not overwhelmed by the precise choice of a prior. It does not seem that the
particular $P(c \mid |Z|)$ adopted in NSB is necessary to have this data-sensitivity, at least as far as the likelihood of $|Z|$
is concerned and at least for the regime tested here.

\subsection{Uncertainty in the Event Space Size}
\label{sec:uncertain_|Z|}

One of the core distinguishing features of the approach analyzed in this
paper is to treat $|Z|$ as a random variable. In particular, the experimental
comparisons with NSB discussed in Section~\ref{sec:post_entropy_results} require specification of
a prior, $P(|Z|)$.

There are several different ways of motivating such a prior over $|Z|$. Perhaps the
simplest is to have it be uniform up to some cutoff, far larger than $\SSS(\n)$. 
A somewhat more sophisticated approach would be to assume that 
$P(|Z|)$ is set by a stochastic process that starts with $|Z| = 1$
and then iteratively adds a new ``bin'' to $Z$, the set of already selected
bins, stopping after each iteration with a probability, $1-\gamma$, ending
at some upper cutoff, ${\underline{m}}$, if we get to $|Z| = \underline{m}$. Then, we can write:
\begin{eqnarray}
P(|Z| \mid {\underline{m}}) &\propto& \gamma^{|Z|}
\end{eqnarray}
for all $|Z| \le \underline{m}$, and $P(|Z| \mid {\underline{m}}) = 0$, otherwise.
Alternatively, to avoid the need to specify an upper cutoff, one could use
$P(|Z|) \propto \exp{(-\alpha |Z|)}$ for some hyperparameter, $\alpha$.

Other models are also possible. For example, we might imagine that there
is some upper bound value, $\underline{m}$, that is set \emph{a priori} and
$|Z|$ initially set to $\underline{m}$. Then, elements of $Z$ are iteratively
removed, at random, stopping after each iteration with a probability, $1-\gamma$, ending
at $|Z| = 1$, if we manage to get to $|Z| = 1$. Whereas the first $P(|Z|)$ is a decreasing function
from $|Z| = 1$ up to $|Z| = \underline{m}$;
this alternative $P(|Z|)$ is an increasing function over that range (see Section~\ref{sec:subset_sel}
for a discussion of these kinds of scenarios where $Z$ is determined by randomly
forming a subset of some original larger set).

A natural extension to both of these models is to allow $\underline{m}$ to vary and
use an associated (hyper)prior.
As always, the choice of random variables and associated (hyper)priors should match
one's understanding of the underlying physical process by which the data is
collected as accurately as~possible.

Once we have a hyperprior, $P(|Z|)$, we can combine it with the
likelihoods plotted in Figure~\ref{est} to get a posterior distribution over $|Z|$. For a $P(|Z|)$
that is nowhere-increasing, since the likelihoods are decreasing functions of $|Z|$, the 
posterior is also a decreasing function of $|Z|$. This would lead to an MAP
 estimate
of $|Z| = |\SSS(\n)|$, though, in general, $\E(|Z| \mid \n)$ would be bigger than $|\SSS(\n)|$.
(In the scenario considered below in Section~\ref{sec:subset_sel}, both the MAP and posterior
expected values of $|Z|$ can exceed $|\SSS(\n)|$.)

\subsection{Specifying a Single Event Space Size} 
\label{sec:single_event_size}

Another way to address uncertainty in $|Z|$ is to set it to a single ``optimal'
value, rather than averaging over it. For example, we could use a 
``coverage estimator''~\cite{bunge1993estimating,vu2007coverage} to set a single value of $|Z|$, $m$ and, then,
evaluate $\E(Q \mid \n, c, |Z| = m)$ using the formulas in WW~\cite{davidfive}. Alternatively, if one has a distribution
over $c$, then we would integrate over it while keeping $|Z|$ fixed. If we assume IUV, so $P(c \mid |Z|) = P(c)$,
this would~give:
\begin{eqnarray}
\E(Q \mid \n, |Z| = m) &=& \int d\rho \; Q(\rho) P(\rho \mid \n, |Z| = m) \nonumber \\
&=& \int d\rho \; Q(\rho) \int dc \; P(\rho, c \mid \n, |Z| = m) \nonumber \\
&=& \int d\rho \; Q(\rho) \frac{\int dc \; P(\n \mid \rho, c, |Z| = m) P(\rho \mid c, |Z| =m) P(c)} {\int dc d\rho \; P(\n \mid \rho, c, |Z| = m) P(\rho \mid c, |Z| =m) P(c)} \nonumber \\
&=& \frac{\int dc d\rho \; P(c) Q(\rho) \big[ \prod_z \rho_z^{n_z -1+ c / m} \big/ G(0, c, m) \big]} 
 {\int dc d\rho \; P(c) \big[ \prod_z \rho_z^{n_z -1+ c / m} \big/ G(0, c, m) \big]}
 \nonumber \\
&=& \frac{\int dc d\rho \; P(c) \Gamma( N + c) Q(\rho) \big[ \prod_z \rho_z^{n_z -1+ c / m)} \big/ [\Gamma(c/m)]^m \big]} 
 {\int dc d\rho \; P(c) \Gamma( N + c) \big[ \prod_z \rho_z^{n_z -1+ c / m} \big/ [\Gamma(c/m)]^m \big]} 
%
\end{eqnarray}
where we have used the results in Section ~\ref{sec:prelim} to derive the last line.
Of course, we could also replace $P(c)$ in this formula with a distribution $P(c \mid |Z| =m)$
if we wish to violate IUV, e.g., as in the NSB~estimator.


\section{Posterior Expected Entropy When $|Z|$ is a Random Variable}
\vspace{-12pt}
\subsection{Fixed $c$}

How can we estimate the entropy of a system when the number of bins is unknown while $c$ is fixed?
Under IUV:
\begin{eqnarray}
P(H = h \mid \vec{n}, c) &\propto& \sum_{|Z| = 1}^\infty \int d \rho \; P(\vec{n} \mid \rho, c, |Z|) 
 P(\rho \mid Z, c) P(|Z| \mid c)~\delta(H[\rho]-h) \nonumber \\
&=& \sum_{|Z| = 1}^\infty \int d \rho \; P(\vec{n} \mid \rho, |Z|)~\D_{c,Z}(\rho) P(|Z|) ~\delta(H[\rho]-h) P(|Z|) \nonumber \\
&\equiv& \sum_{|Z| = 1}^\infty \int d \rho \; P(\vec{n} \mid \rho)~\D_{c,Z}(\rho) P(|Z|) ~\delta(H[\rho]-h)
\end{eqnarray}
where $P(|Z|)$ could be given by one of the priors discussed in Section~\ref{sec:uncertain_|Z|}.

Instead of trying to solve for this, as in WW, we can consider the posterior expected
values of the moments of the entropy. Using IUV, the first moment is:
\begin{eqnarray}
\label{avgent}
\E(H \mid \vec{n}, c) &=& \sum_{|Z| = 1}^\infty\int d\rho~ H(\rho) P(\rho, |Z| \mid \vec{n}, c) \nonumber \\
&=& \sum_{|Z| = 1}^\infty \int d\rho~ H(\rho) \frac{ P(\vec{n} \mid \rho) P(\rho \mid |Z|, c) P(|Z|) }
{\sum_{|Z'| = 1}^\infty \int d\rho' ~ P(\vec{n} \mid \rho') P(\rho' \mid |Z'|, c) P(|Z'|)} \nonumber \\
&=& \frac{\sum_{{|Z|}=M}^\infty P(|Z|) \int d\rho~H(\rho)~P(\vec{n} \mid \rho)~\D_{c,Z}(\rho)}
{\sum_{{|Z|}=M}^\infty P(|Z|) \int d\rho ~ P(\vec{n} \mid \rho)~\D_{c,Z}(\rho)}
\end{eqnarray}

The two integrals are straight-forward. We already know the denominator; the numerator is not that much harder. 
To simplify the expression of the result, define $M \equiv |\SSS(\n)|$ and write $\{n_i : i \in \SSS(\n)\}$ for the set, $\{n_Z(z) : z \in \SSS(\n)\}$.
Then, we get:
{\small
\begin{eqnarray}
\E(H \mid \vec{n}, c) &=& \frac{1} {\sum_{{|Z|}={M}}^\infty P(|Z|) \frac{\prod_{i=1}^{{M}} \Gamma(n_i+c/{|Z|})}{\Gamma(c/{|Z|})^M}} \; \times \nonumber \\
&& \:\:\:\: \bigg\{ \sum_{{|Z|}=M}^\infty P(|Z|) \frac{\prod_{i=1}^{M} \Gamma(n_i+c/{|Z|})}{\Gamma(c/{|Z|})^M} 
\bigg(\bigg[\sum_{i=1}^{M} \frac{n_i+c/{|Z|}}{N+c} \Delta\Phi^{(1)}(n_i+c/{|Z|}+1,c+N+1)\bigg] \nonumber \\
&& \:\:\:\:\:\: \:\:\:\:\:\: \:\:\:\:\:\:\:\:\:\:\:\: \:\:\:\:\:\: \:\:\:\:\:\:\:\:\:\:\:\:\:\:\:\:\:\:\:\:\:\:\:\: +~ \:\: ({|Z|}-M) \frac{c/{|Z|}} {N+c} 
\Delta\Phi^{(1)}(c/{|Z|}+1,c+N+1)
\bigg) \bigg\} \nonumber \\
\label{eq:post_ent}
\end{eqnarray}}
where the term on the last line arises from empty bins and where $\Phi^{(n)}$ and $\Delta\Phi^{(n)}$ are defined in Section~\ref{sec:prelim}. 

Recall from Section~\ref{sec:|Z|_like} that the likelihood, $P(\n \mid c, |Z|)$, is
strongly weighted towards small $|Z|$. Therefore,~if $P(|Z|)$ has an upper
cutoff, $\mathscr{M}$, and is nowhere-increasing, then under a DI hyperprior, the posterior, 
$P(|Z| \mid \n, c)$, must also be strongly weighted towards small $|Z|$. This will cause
posterior moments of $\rho$ to be dominated by values, $|Z|$, that are not much
larger than $M$. In general, this will mean that these posterior moments will not
be prior-dominated, in the sense that attributes of $\n$ will affect them significantly.

As an illustration, note that, in general, the innermost sum in Equation~(\ref{eq:post_ent}):
\begin{eqnarray}
\bigg[\sum_{i=1}^{M} \frac{n_i+c/{|Z|}}{N+c} \Delta\Phi^{(1)}(n_i+c/{|Z|}+1,c+N+1)\bigg] ~+~ ({|Z|}-M) \frac{c/{|Z|}} {N+c} \Delta\Phi^{(1)}(c/{|Z|}+1,c+N+1) \nonumber
\end{eqnarray}
goes to an $\n$-dependent constant as $|Z|$ becomes large.
(The arguments of both $\Delta \Phi^{(1)}$'s become independent of $|Z|$, and the factors
multiplying each of them go to a constant.) Furthermore, as discussed in Section~\ref{sec:|Z|_like}, $\Gamma(c/|Z|)^{-M}$ 
goes to zero as $|Z|$ gets large. Accordingly, for our assumed form of $P(|Z|)$, once $\mathscr{M}$ is appreciably larger than $M$,
posterior expected entropy does not change if $\mathscr{M}$ instead becomes hugely
larger than $M$. In this sense, the posterior expected
entropy does not become prior dominated as $\mathscr{M}$ becomes very large.
(At a minimum, $M$---an attribute of the data, not the prior---is determining
the range of relevant $|Z|$.) Therefore, to the degree that prior-dominance 
is avoided in Equation~(\ref{eq:post_ent}), treating $|Z|$ as a random variable and requiring IUV
removes the phenomenon that caused NSB to adopt their scheme for setting $P(c)$.
%

A formula similar to Equation~(\ref{eq:post_ent}) gives the second moment of the posterior distribution over the entropy. 
(It is too long to write out here; we recommend that a package like Mathematica be used to evaluate it.) Combining 
that formula with Equation~(\ref{eq:post_ent}) provides the posterior variance of the entropy when
the number of bins is a random variable.

\subsection{Uncertain $c$}
\label{sec:post_entropy_results}

To allow for varying $c$, one can change the sum over ${|Z|}$ in Equation~(\ref{eq:post_ent})
for a fixed $c$ to a new expression involving sums over ${|Z|}$ together with integrals over $c$ with some prior for $c$. Care
must be taken when doing this, since $c$ is not independent of $\n$ (see Section~\ref{sec:single_event_size} for another example
of integrating over $c$ when conditioning on $\n$). Using IUV, the result is:
\begin{eqnarray}
\E(H \mid \vec{n} ) &=& \sum_{|Z| = 1}^\infty\int d\rho dc ~ H(\rho) P(\rho, |Z|, c \mid \vec{n}) \nonumber \\
&=& \frac{\sum_{{|Z|}=M}^\infty P(|Z|) \int d\rho dc~H(\rho)~P(\vec{n} \mid \rho)~\D_{c,Z}(\rho) P(c)} 
{\sum_{{|Z|}=M}^\infty P(|Z|) P(c) \int d\rho dc ~ P(\vec{n} \mid \rho)~\D_{c,Z}(\rho) P(c)}
\label{awesome-equation}
\end{eqnarray}
To evaluate this, we need only separately apply $\int dc P(c)$ to both the numerator and denominator in Equation~\eqref{eq:post_ent}.

\subsection{Experimental Tests}

The primary focus of this paper is an analysis of the IUV desideratum and its implications
for the hyperprior. However, as a sanity check, in this subsection, we compare the performance of posterior
estimators of entropy and of mutual information that are based on a DI hyperprior to three estimators 
of those quantities that were previously considered in the literature.
These three alternative estimators are NSB, the estimator considered in~\cite{nemenman2011coincidences}
(which is an asymptotic version of NSB that 
allows for the estimation of entropy when the number of bins is unknown), and the ``Coverage-Adjusted Estimator'' (CAE) of~\cite{vu2007coverage}. To simplify the exposition, we will sometimes refer to any estimator based on a
DI hyperprior as a ``W\&D'' estimator.

From a decision-theoretic perspective, no estimator can do better under an assumed
hyperprior than one that is Bayes-optimal for that hyperprior, \emph{if} one quantifies performance
with experiments that draw samples from that same hyperprior.
Failures of a Bayes-optimal estimator always arise from a mismatch between the hyperprior used to construct the estimator and the hyperprior used to actually generate the data.
Accordingly, to make meaningful comparisons between a W\&D estimator
and others, we must see how well Equation~(\ref{awesome-equation}) performs ``out of class'', for
data that is not generated from the DI hyperprior used to construct the W\&D estimator.

To make these comparisons, we use a
W\&D estimator for a hyperprior that has logarithmic $P(c)$ and uniform $P(|Z|)$
and then consider two general types of data
that are \emph{not} drawn under this hyperprior. In particular,
we consider: (1) distributions sampled from Dirichlet distributions with fixed $c$; and (2) power-law distributions of the form:
\begin{equation}
P(i) \propto \frac{1}{S[i]^\alpha},~[i=1\ldots m]
\label{eq:power_law}
\end{equation}
where $S[i]$ is a one-to-one map on the integers from one to $m = |Z|$.
(Varying such $S[.]$ ensures that the order of the terms is not fixed, but can vary depending on the particular choice of $S[.]$; this is important when constructing joint probability distributions by re-interpreting a probability distribution over $m$ categories as a joint distribution over $\sqrt{m}\times\sqrt{m}$ categories, as in estimates of mutual information between a pair of random variables.)

We estimated the entropy and mutual information based on datasets
generated this way using Equation~\eqref{awesome-equation}, the Coverage-Adjusted Estimator of~\cite{vu2007coverage} (Equation~18, with $n + 1$ correction to make it well-defined in the singleton case), the Asymptotic NSB estimator of~\cite{nemenman2011coincidences} and a ``large-$Z$'' version of the standard NSB estimator of~\cite{Nemenman:2002fk}, where the bin size is set to a large value assumed to be larger than the possible number of bins. For both W\&D and the large-$Z$ NSB, we must include a maximum bin number. We take this to be $10,000$ ({\em i.e.}, a hundred times larger than the actual number), and in the case of mutual information estimation, we allow the large-$Z$ NSB estimator to assume that the maximum bin number for each marginal is 100 ({\em i.e.}, ten times larger than the actual one), with the maximum bin number for the full distribution, as before, set to $10,000$.

We then computed the {\sc RMS
} error between these estimates and the truth for all three estimators. 
The results are shown in Figure~\ref{fig:post_entropy}, where we sample 100-bin processes in the ``deeply undersampled regime'', where $N$, the number of counts, is $\sqrt{m}$.

\begin{figure}[H]
\centering
\includegraphics[width=4in]{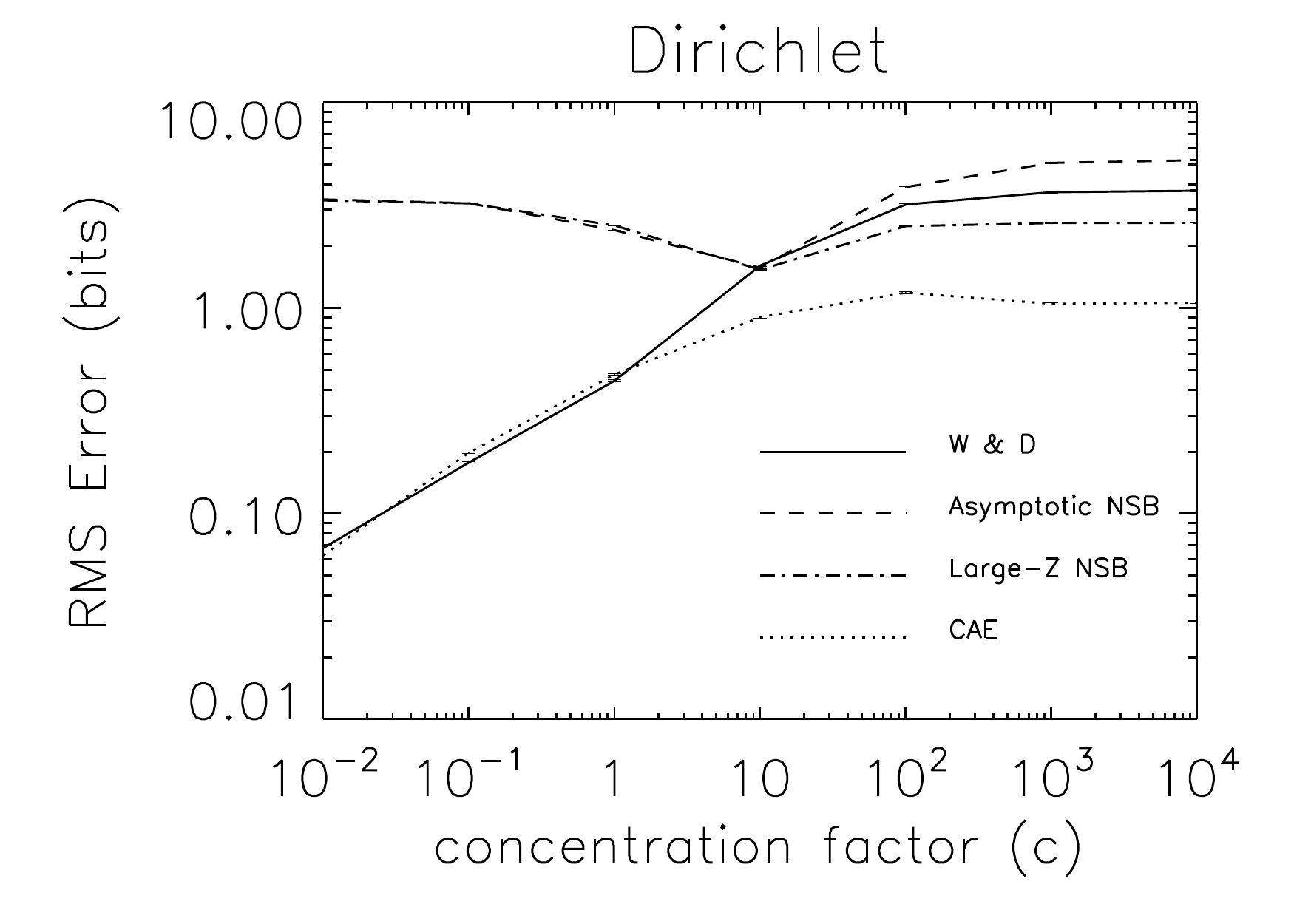} 
\caption{RMS 
 error for the estimation of entropy with an unknown bin number. 
In both plots, the solid line shows the W\&D estimator, Equation~(\ref{awesome-equation}); the dashed line 
shows the asymptotic version of NSB 
 (Equation~(29) of~\cite{nemenman2011coincidences}); the dotted line 
shows the Coverage-Adjusted Estimator of~\cite{vu2007coverage}; the dot-dashed line shows the ``large-$Z$'' version of NSB.
The true bin number is 100, and the associated
distribution, $\rho$, was sampled ten times ({\em i.e.}, it was radically under-sampled given the number of bins). ({\bf Top}) 
Results when the distribution, $\rho$, used to generate the data is
randomly generated under a Dirichlet distribution, whose concentration, $c$, is varied from $10^{-2}$ (highly non-uniform; many hits in a few bins) to $10^4$ (highly uniform), as indicated. ({\bf Bottom}) Results when $\rho$ is a
a power-law distribution over the bin number with index $\alpha$ that is varied from zero (perfectly uniform) to four (highly non-uniform),
as in Equation~\eqref{eq:power_law}. (The Zipf distribution corresponds to $\alpha = 1$.) }
\label{fig:post_entropy}
\end{figure}

\begin{figure}[H]\ContinuedFloat
\centering	
\includegraphics[width=4in]{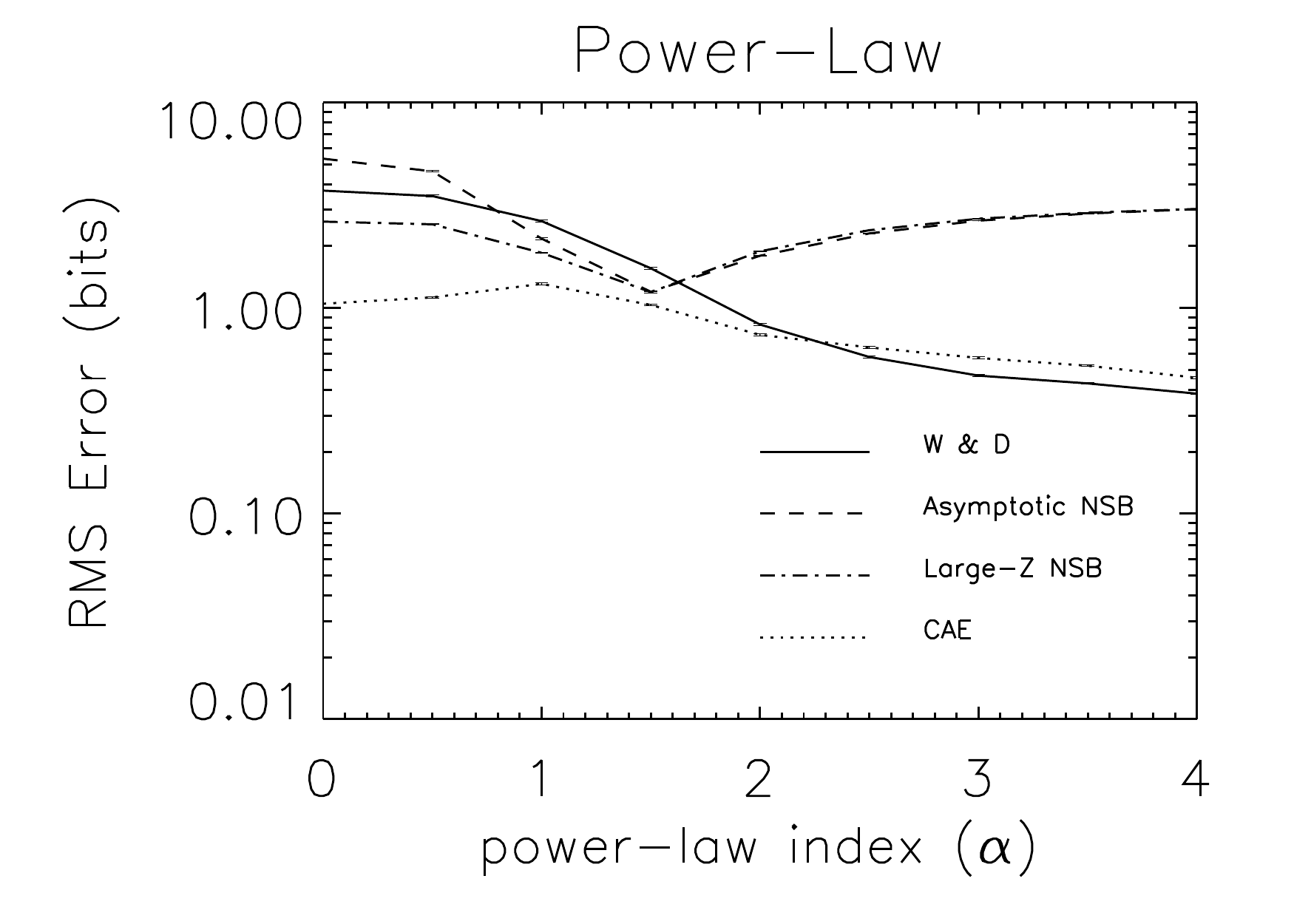}
\caption{\em Cont.}
\end{figure}

The estimator of Equation~(\ref{awesome-equation}) outperforms the asymptotic NSB estimator for a wide range of distributions, both Dirichlet and power-law. This is not entirely unexpected; the asymptotic estimator works only to zeroth order in $1/N$ and $1/m$. In the regimes
where the true $\rho$ is likely to have low-entropy, the large-$Z$ estimator performs almost identically to the asymptotic NSB estimator, while it performs somewhat better in the high-entropy regime, where it is competitive with Equation~(\ref{awesome-equation}). For low-entropy samples (either low $c$ or high $\alpha$), Equation~(\ref{awesome-equation}) is competitive with the Coverage-Adjusted Estimator. Inefficiencies in Equation~(\ref{awesome-equation}) trace back to our prior over $c$; in cases where one has a strong belief that the data are drawn from high-entropy distributions, different $c$ weights should be used.

Interestingly, at the $\alpha$ equal to the unity point (the Zipf distribution), all four methods are within a factor of two of each other {\sc rms
}. The strongest differences between the methods emerge at low entropies.

%

\begin{figure}[H]\centering
\includegraphics[width=4in]{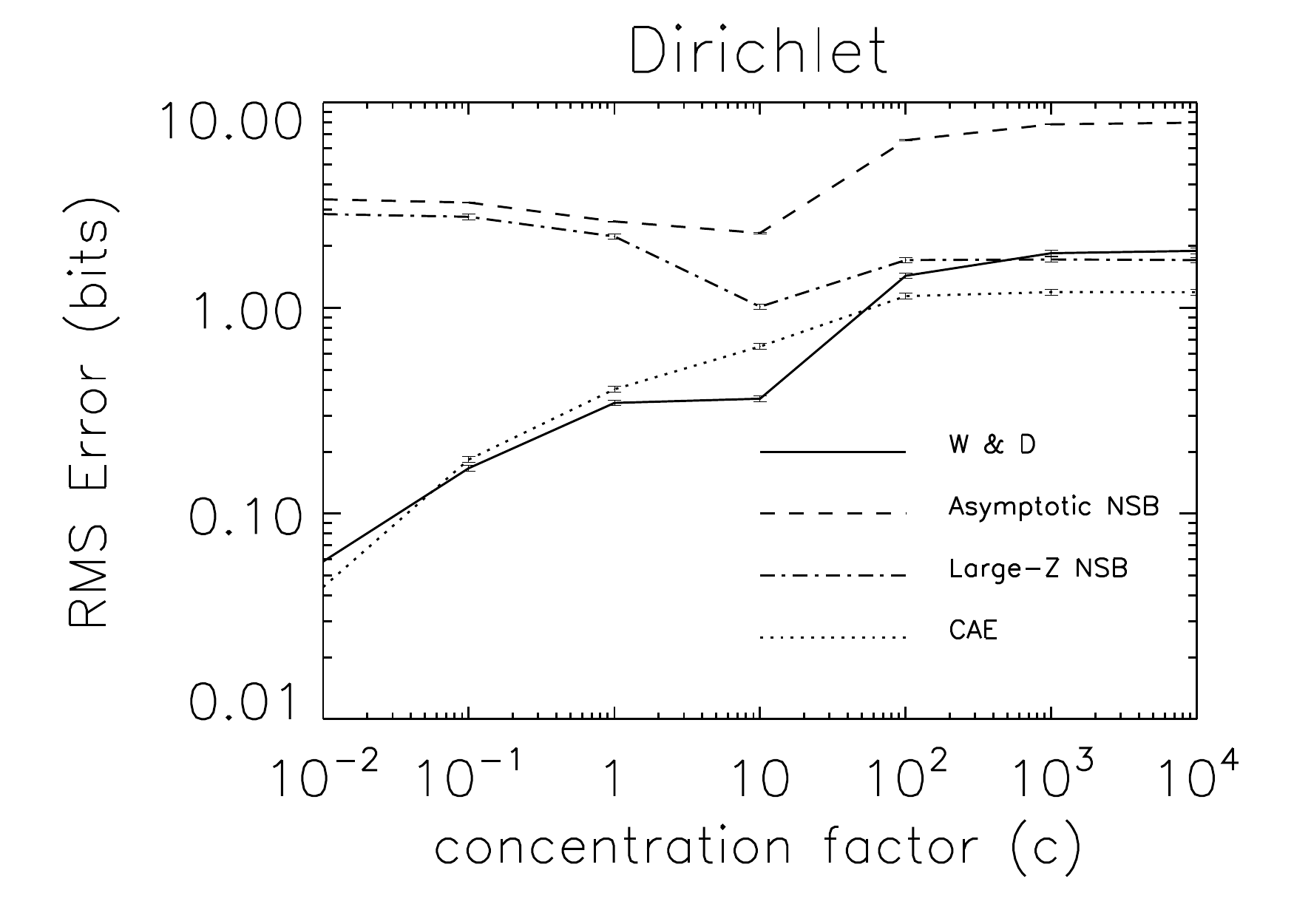} \\
\includegraphics[width=4in]{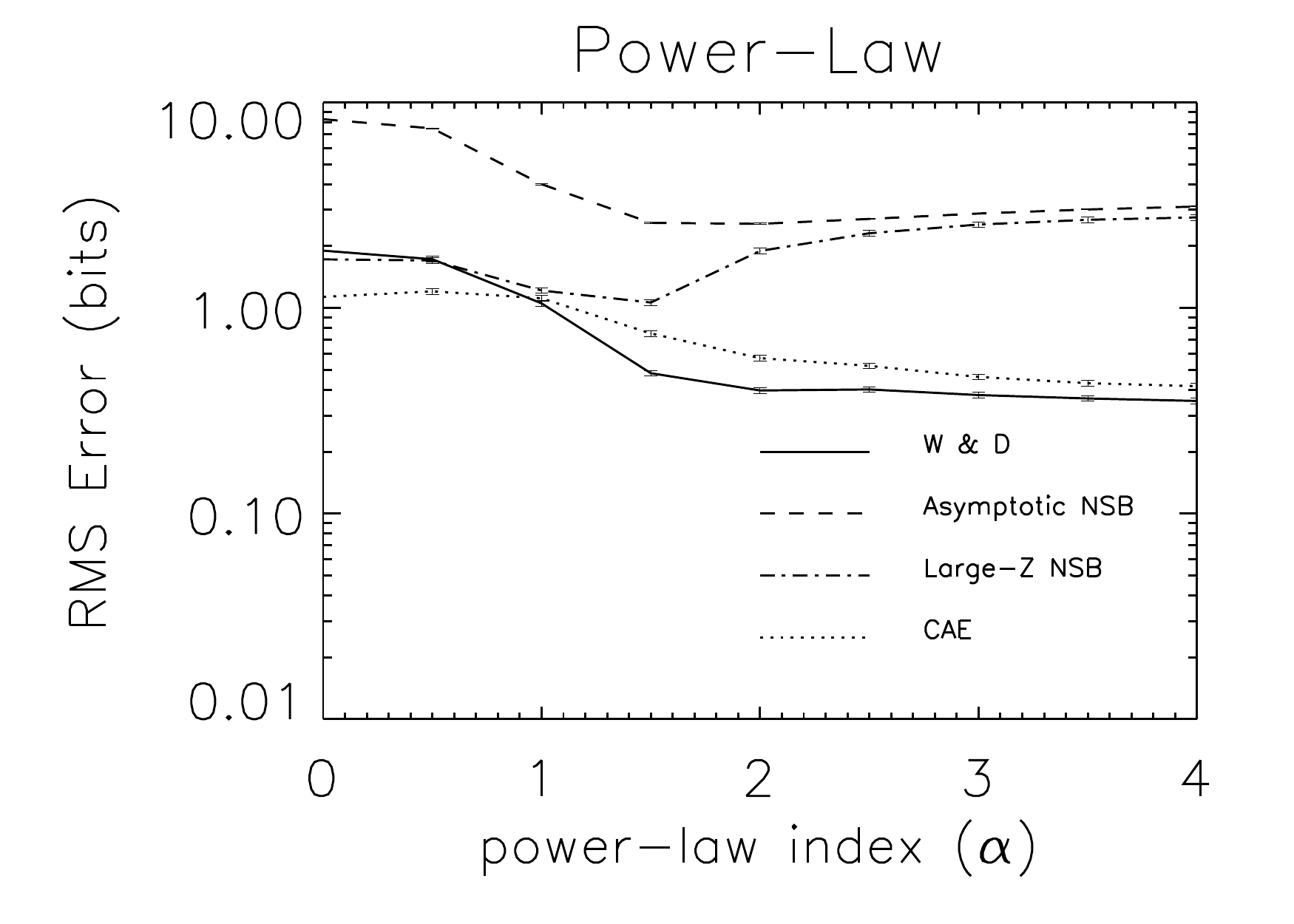}
\caption{RMS
 error for the estimation of mutual information with unknown bin number. The true bin number is 100, interpreted as events over a $10\times10$ joint space, sampled ten times ({\em i.e.}, radically under-sampled compared to the number of bins).
The solid line shows the estimate made using the W\&D estimator, Equation~(\ref{awesome-equation}); the dashed line, the asymptotic version of NSB (Equation~(29) of~\cite{nemenman2011coincidences}); the dotted line, the Coverage-Adjusted Estimator of~\cite{vu2007coverage}; the dot-dashed line, the ``large-$Z$'' version of NSB. ({\bf Top}) Results when $\rho$ is
randomly generated under a Dirichlet distribution whose concentration, $c$, is varied from $10^{-2}$ (highly non-uniform; many hits in a few bins) to $10^4$ (highly uniform), as indicated. ({\bf Bottom}) Results when $\rho$ is
is a power-law distribution over the bin number with index $\alpha$, which is varied from zero (perfectly uniform) to four (highly non-uniform). (The Zipf distribution corresponds to $\alpha = 1$.)}
\label{fig:post_mutual_info}
\end{figure}

We can use the same methods to compare the accuracies of the estimators of the mutual information. In particular, since Equation~(\ref{awesome-equation}) respects IUV, we can decompose the mutual information into the sum and differences of entropies. In Figure~\ref{fig:post_mutual_info}, we plot {\sc RMS} error for mutual information estimated in this fashion and compare it to the naive use of the asymptotic and large-$Z$ NSB estimators and the Coverage-Adjusted Estimator of~\cite{vu2007coverage}.

The differences between the estimators of mutual information are more extreme
than the difference for the estimators of entropy; the W\&D estimator based on Equation~(\ref{awesome-equation}) and the Coverage-Adjusted Estimator perform comparably. The large-$Z$ NSB Estimator performs comparably in the high-entropy regime; the asymptotic NSB estimator tends to perform poorly.

We emphasize
that the choice of DI hyperprior used in these experiments was ``naive'', not based on any careful reasoning or
desiderata. In particular, no consideration was given to what type of generative process 
plausibly underlies the construction of $|Z|$ (see Section~\ref{sec:generative}).
Potentially more compelling results would occur for a
more careful choice of hyperprior.

\section{Generative Models of $Z$}
\label{sec:generative}

In this section, we discuss subtleties in how one models the statistical generation of $Z$.
To ground the discussion, we will sometimes consider an example where
a fishery's biologist randomly samples fish from a lake via a catch-and-release
protocol, to try to ascertain quantities, like the number of fish species in the lake,
the entropy of the distribution of the counts of members of those species, \textit{etc}. In this example, $Z$ is the 
set of all fish species in the lake (and is explicitly a random variable),
and $\vec{n}$ is a set of the counts of the species of fish.

\subsection{Mapping Physical Samples to Bin Labels}

In the processes considered in the previous sections, $\pi(m)$ is sampled to produce
$m$, and a set, $Z$, of $m$ elements, labeled, for example, $\{1\ldots, m\}$, is
then created. Next, $c$ is sampled from $\pi(c \mid m)$. After this,
$P(\rho \mid c, m)$ is sampled to get a $\rho$, and finally, $\vec{n}$ is
sampled from $\rho$. 

Note that the $\vec{n}$ produced at the end of this process is a set of
counts of the integers ranging from one to $m$. However, physically, $\n$
is not a set of counts of integers. This implies that we need a map
from the physical characteristics of the samples in the real world into $\{1, \ldots, m\}$. As an illustration, \linebreak in the fish-in-a-lake example,
$\vec{n}$ is a set of the counts of species of fish that are distinguished by
their physical characteristics (assuming no DNA sequencing or the like is used). Therefore, to
apply the formulas derived in the previous sections, the biologist provided with a sample of the counts of fish species needs an
invertible map sending each species of fish in their sample into $\{1\ldots, m\}$.

How should the biologist create that map? One idea might
be to randomly build an invertible map taking each of the distinct
species of fish they have sampled to a different member of the set of integers from one to $m$. However,
the biologist cannot do this, since they do not know $m$, and, so, cannot build such a map. ($m$ is a random variable,
whose value is not known with certainty to the biologist, even after the biologist gets $\n$). Another possibility
would be to assign the species of the first fish the biologist samples to one, the
second distinct species sampled to two, and so on. However, this would
introduce major biases in the estimators (for example, before \emph{any} data is generated,
we would know that $n_1 \ge 1$, something we would not know for any $n_i$, where $i > 1$).

%
As an alternative,
we can model the statistical
process as one in which the biologist 
%
%
assigns a species ``label'' to each new fish as the biologist draws it from the lake, based on the physical
characteristics of that fish. More precisely, assuming the biologist measures
$K$ real-valued physical characteristics of each fish that the biologist draws from the lake, we can model
the sampling process as follows:

\begin{enumerate}

\item $\pi(m)$ is sampled to get $m$, the number of fish species in the lake. At this
 point, nothing is specified about the physical characteristics of each of those $m$ species.

\item Next, a Dirichlet distribution is sampled that extends over distributions $\rho$ that themselves are defined
 over $m$ bins.

\item Next, a vector, $v_j \in R^K$, is randomly assigned to each of $j = 1, \ldots, m$,
 e.g., where each of the $m$ vectors is drawn from a Gaussian
 centered at zero (the precise distribution does not matter).
 $v_j$ is the set of $K$ real-valued physical characteristics that we
 will use to define an idealized canonical specimen of fish species, $j$.
By identifying the subscript, $j$, on each $v_j$ as the associated bin integer,
we can view $\rho$ as defined over
 $m$ separate $K$-dimensional vectors of fish species characteristics.

\item $\rho$ is IID sampled to get a dataset of counts for each species, one through $m$.
 Physically, this means that the biologist draws a fish from bin $j$ with probability
 $\rho_j$, {\em i.e.}, they draw a fish with characteristics $v_j$ with probability $\rho_j$.
 
\end{enumerate}

Note that we could interchange the order of steps 2 and 3.
Note also that, in practice, since lakes do not contain ``ideal fish'', there will be some small
 noise added to $v_j$ each time the biologist draws a member of species, $j$. We
can assume that noise is small enough on the scale of the typical distance between the vectors of
 canonical fish species characteristics, so that the probability is infinitesimal,
 that the biologist misassigns what species a fish the biologist draws belongs to. 

This generative model is more elaborate than the more informal one described at the beginning of this
section. However, both models result in the same formulas, namely, those given in the previous sections.

\subsection{Subset Selection Effects}
\label{sec:subset_sel}

There are other, simpler models that one might think solve the difficulty of how to map the physical
members of $\n$ into integers $\{1, \ldots, \m\}$. However, many of these alternative models
introduce subtle biases into the estimators. Seemingly trivial differences 
in the formulation of the problem of estimating functionals from data can have substantial
effects on the resultant predictions. 

To illustrate this, consider the common scenario where
we are given a set, $\hat{Z}$, and know that $Z \subseteq \hat{Z}$, but do not
know which precise subset of $\hat{Z}$ is $Z$. A simple example of such a scenario
is a variant of the one where a field biologist wishes to estimate the entropy of the fish
species in a particular lake by IID sampling fish in that lake. Say the biologist knows the set
of all fish on Earth. However, assume they have no \emph{a priori} knowledge
that one fish is more likely than another to be a lake-dwelling fish. 
In this case, $\hat{Z}$ is the subset of all fish on Earth,
and $Z$ is the set of all lake-dwelling fish. While the biologist knows $\hat{Z}$,
they are uncertain of $Z$. 

We might presume that the estimation of quantities, like posterior expected
entropy of the distribution of fish in the lake, would not depend on whether we
calculate them for this scenario where we know that $Z$ is an (unknown)
subset of $\hat{Z}$ or, instead, calculate them for the original scenario analyzed above,
where we simply have uncertainty of $Z$, without any concern for an embedding
set, $\hat{Z}$. However, it turns out that
the estimation is quite different in these two scenarios. This illustrates how
much care one must take in the statistical formulation of the estimation problem.

To see that estimation differs in this variant of the fish-in-a-lake example, first, 
as shorthand, define $m$ to be the size of $Z$, $|Z|$. In the subset-of-$\hat{Z}$ scenario:
\begin{eqnarray}
P(\n \mid \hat{Z}) &=& \sum_{m = |\SSS(\n)|}^{|\hat{Z}|} P(m \mid \hat{Z}) P(\n \mid m, \hat{Z}) 
\label{eq:new_like}
\end{eqnarray}
(Note the implicit definition of $\n$ as the vector of counts for all elements in $ \hat{Z}$.)
$P(m \mid \hat{Z})$ plays the same role here as the prior, $P(|Z|)$, does in the analysis above of the
original scenario where there is no $\hat{Z}$. 

To evaluate the likelihood $P(\n \mid m, \hat{Z})$ in Equation~\eqref{eq:new_like}, we need a stochastic model of 
how $Z$ is formed from $\hat{Z}$. There are many such models possible. For simplicity, 
adopt the model that all $Z$'s of a given size, $m$, are equally likely:
\begin{eqnarray}
P(Z \mid m, \hat{Z}) &=& \frac{\delta_{|Z|, m}} {\binom{|\hat{Z}|} {m}} 
\label{eq:subset_sel_model}
\end{eqnarray}
Recalling the definitions of $G$ and $\SSS$ in Section~\ref{sec:prelim}, we have the following
(the proof is in Appendix B):
\begin{proposition}
Under the conditional distribution in Equation~\eqref{eq:subset_sel_model}:
\begin{enumerate}
\item $P(\n \mid Z, m, \hat{Z}) \; \propto \; I(\SSS(\n) \subseteq Z) G(\n, c, m)$ ;
\item $P(\n \mid m, \hat{Z})\; \propto \; \binom{|\hat{Z}| - |\SSS(\n)|}{m - |\SSS(\n)|} \; G(\n, c, m) $
\end{enumerate}
\label{prop:new_like}
\end{proposition}

In contrast to the likelihood $P(\n \mid m, \hat{Z})$ given by Proposition~\ref{prop:new_like}
for the subset-selection scenario, the likelihood for the original scenario analyzed above was:
\begin{eqnarray}
P(\n \mid m) &=& \int d\rho_Z \; P(\n \mid \rho_Z) P(\rho_Z \mid Z) \nonumber \\
&\propto& G(\n, c, m) 
\end{eqnarray}
where $|Z| = m$.
Therefore, writing them out in full:
\begin{eqnarray}
P(\n \mid m) &=& \frac{G(\n, c, m) }{\sum_{\n'} G(\n', c, m)}
\end{eqnarray}
\begin{eqnarray}
P(\n \mid m, \hat{Z}) &=& \frac{ \binom{|\hat{Z}| - |\SSS(\n)|}{m - |\SSS(\n)|} \; G(\n, c, m) }{\sum_{\n'} \binom{|\hat{Z}| - |\SSS(\n')|}{m - |\SSS(\n')|} \; G(\n', c, m) }
\end{eqnarray}
(The sum in the denominator in the second equation is implicitly restricted to those
$\n' \in \hat{Z}$ whose support contains no more than $m$ elements, and in both sums,
we are implicitly restricting attention to those $\n'$ with the same total number of elements as $\n$.)
Intuitively, the reason for the difference between these two likelihoods
is that in the subset-of-$\hat{Z}$ scenario, there are combinatorial effects reflecting
the number of ways of assigning elements in $\hat{Z} \setminus \SSS(\n)$ to the $|Z| - |\SSS(\n)|$ bins
in $Z$ that are unoccupied in $\n$, whereas there are no such effects in the original scenario.

The extra combinatoric factor, $\binom{|\hat{Z}| - |\SSS(\n)|}{m - |\SSS(\n)|}$, in the subset-of-$\hat{Z}$ scenario's likelihood 
for $m$ pushes that likelihood
to prefer smaller values of $|\SSS(\n)|$ compared to the original likelihood. It also distorts
how the likelihood depends on $m$. To a degree, we can compensate for this second effect using
the other term 
in the summand in Equation~\eqref{eq:new_like} besides $P(\n \mid m, \hat{Z})$, $P(m \mid \hat{Z})$.
Even once we do this though, the~precise estimates generated in the two scenarios will differ in general,
since the prior, $P(m \mid \hat{Z})$, cannot fully compensate for an effect that depends on the data, $\n$.

More generally, one might even want to treat $|\hat{Z}|$ as a random variable. Returning
to our concrete example involving fish, we would do this if we are uncertain about the total
number of fish on Earth. This re-emphasizes the point that,
as always, the choice of random variables and associated priors should match
one's understanding of the underlying physical process by which the data is
collected as accurately as possible.

%
%
%

\section{Conclusions}

The problem of estimating a functional of a distribution, $\rho$,
based on samples of $\rho$ is a core concern of statistics.
In particular, in recent decades, there has been a great deal of work on estimating \linebreak information-theoretic functionals
of $\rho$ based on samples of $\rho$.

Bayesian approaches to this problem began with WW, where the Dirichlet prior for
$\rho$ was adopted. This work concentrated on the case where the concentration
parameter, $c$, of the Dirichlet prior equals the size of the underlying
event space, $|Z|$. NSB pointed out that this special case has the problem 
that the resultant estimators are prior dominated whenever $|Z|$ is much larger
than the support of the dataset. NSB realized that this problem could
be addressed by using a hyperprior over $c$. They then advocated an
approach to setting $P(c)$. Unfortunately, there is a substantial
problem with the choice $c = |Z|$ analyzed in WW that is not
fixed by the NSB choice of $P(c)$. In both approaches, the posterior expected
value of a quantity, like mutual information, will depend on which of
the many equivalent definitions of mutual information one adopts.
In other words, that posterior expected value of such quantities is ill-defined.

In many situations, there will be 
uncertainty of $|Z|$, as well as $c$. Indeed, arguably, 
there is \emph{always} an uncertain
number of hidden degrees of freedom in the stochastic process
that produced the data, degrees of freedom not recorded as
components of that data. Since the stochastic process model must be set independently
of the likelihood and it is the likelihood that determines
what degrees of freedom are recorded, we must allow $\rho$
to run over those hidden degrees of freedom, as well as the
visible ones recorded in the data. Since we typically do not know how many such hidden
degrees of freedom there are, this means we have an uncertain
value of $|Z|$.

This reasoning argues that we should use a hyperprior, $P(c, |Z|)$. 
To do so, we must specify how the concentration parameter
is statistically coupled with the size of the underlying event
space. It is not at all clear how to do that in a hierarchical Bayesian way, where
we cannot consider either the likelihood (which determines what variables
are observed) or how the posterior estimate of $\rho$ would be used (which is what
NSB uses to couple $c$ and $|Z|$). It is also not at all
clear how to specify a prior that extends over hidden degrees of freedom.

In this paper, we address the second concern by introducing the desideratum that
\emph{for any functional that only depends on those components of $\rho$
corresponding to the recorded degrees of freedom}, the number of hidden degrees of freedom has no effect
on our estimate of the functional. This desideratum says that
our second problem is not a problem. We prove that this ``Irrelevance of 
Unseen Variables'' (IUV) desideratum can be satisfied, but only if $c$ and $|Z|$ are independent. Therefore, IUV resolves both of our~concerns. 

In deriving this result, we prove an intermediate result that simplifies the calculation of
some posterior moments. In particular, we show how to use it to derive the formula for
posterior expected mutual information given in WW in essentially a single line.

We also show that by using a $P(c, |Z|)$ consistent with IUV rather than
the one used in either WW or NSB, we resolve the problem shared by them
that posterior expected mutual information is ill-defined. In addition,
as we illustrate, using a hyperprior that respects IUV can also greatly simplify calculation of posterior moments
of information-theoretic functionals. Another advantage of allowing
$|Z|$ to vary and adopting IUV's hyperprior is that
posterior expected values of information-theoretic quantities are no longer 
prior-dominated as they are under the hyperprior of WW. In this sense, there is no need for using
approximations for setting $P(c)$, as, for example, the scheme in NSB.

After presenting these results, we discussed both hierarchical Bayesian approaches and other approaches
for estimating information theoretic quantities when $m$ and $c$ are both random variables.
We ended by discussing some changes to the statistical formulation of the 
estimation problem that would appear to be innocuous, but can actually substantially
affect the resultant estimates.

\section*{\noindent Acknowledgments}
\vspace{12pt}

We would like to thank John Young, Michael Hurley and Gordon Pusch for
their assistance in compiling the errata. S.D. acknowledges the support of the Santa Fe Institute Omidyar Postdoctoral Fellowship, the National Science Foundation Grant EF-1137929, ``The Small Number Limit of Biological Information Processing'' and the Emergent Institutions Project. D.H.W. acknowledges the support of the Santa Fe Institute.

\section*{\noindent Conflicts of Interest} 
\vspace{12pt}

{The authors declare no conflict of interest.}

\appendix

\section*{\noindent Appendix A---Relevant Results and Errata from WW }

$ $

In this appendix, we review relevant results from WW for the case of fixed $Z$ and $c = |Z|$ and present
errata for those results in WW. 
(A preliminary set of errata was reported in~\cite{wowo_erratum}.)

%

It is straight-forward to generalize the reasoning in WW to derive:
\begin{equation}
\E(H \mid \n) = -\sum_z \frac{n_z+c/|Z|}{N+c}\Delta\Phi^{(1)}(n_z+1+c/|Z|, N+c+1)
\label{eq:post_ent_mean}
\end{equation}
where care must be taken to replace the quantity ``$n_i$'' in WW with ``$\n(x, y) - 1 + c / |Z|$'',
since WW considered the uniform (Laplace) prior, in which $c = |Z|$ and $L(.)$ is flat.
Continuing to make the assumption of WW (and many others) that $L$ is flat, we
can evaluate the posterior variance for arbitrary $c$ as:
\begin{eqnarray*}
\E( [H - \E(H \mid \n)]^2 \mid n)
& = & \sum_{z\neq z'}\frac{(n_z+1)(n_{z'}+1)}{(N+c)(N+c+1)} A_{z,z'}\\
& & +\sum_z \frac{(n_z+1)(n_z+2)}{(N+c)(N+c+1)}B_z
\end{eqnarray*}
where:
\begin{eqnarray}
A_{z,z'} & = & \Delta\Phi^{(1)}(n_z+2,N+c+2) \label{phierror} \\
& & \times\Delta\Phi^{(1)}(n_{z'}+2,N+c+2) \nonumber \\
& & -\Phi^{(2)}(N+c+2) \nonumber
\end{eqnarray}
and:
\begin{eqnarray*}
B & = & [\Delta\Phi^{(1)}(n_z+3,N+c+2)]^2 \\
& & +\Delta\Phi^{(2)}(n_z+3,N+c+2) \\
\end{eqnarray*}
where $\Phi^{(n)}$ and $\Delta\Phi^{(n)}$ are defined in Section~\ref{sec:prelim}. 

The variance was incorrectly reported in WW: there was an error in its version of the second line of Equation~\eqref{phierror}.
The mean expected entropy, in the absence of data and under the Laplace prior, is:
\begin{equation}
\E(H) = -\Delta\Phi^{(1)}(2,m+1) = \sum_{q=2}^m \frac{1}{q}
\end{equation}
(This was incorrectly reported in~\cite{wowo_arxiv}, but correct in WW.) 

A complete list of errata in the published article follows. Errata unique to the arXiv versions are not~shown.
\begin{enumerate}
\item The Dirichlet prior equation in the continued paragraph on page 6843 should have the summation symbol replaced with the product symbol.
\item Theorem 8 on page 6846---error as described above, corrected in Equation~\eqref{phierror}. The analogous equation, $\mathcal{E}_{\overline{IJMN}}$ (WW1, page 6852), does not contain the analogous error.
\item Definitions necessary for various subsets, on page 6851, have errors. In particular, $\nu_i$ should be $n_i+1$, and $\gamma_n$
should be $\prod_{i=1}^n\Gamma(\nu_i)$.
\item There is an error in the definition of $\mathcal{E}_{\overline{IN}}$ (page 6852.) In particular, the $\nu$ symbols in the denominators of the term:
\begin{displaymath}
\left(1-\frac{\nu_{i\cdot}+\nu_{\cdot n}-2\nu_{in}}{\nu}+\frac{(\nu_{i\cdot}-\nu_{in})(\nu_{\cdot n}-\nu_{in})}{\nu(\nu+1)}\right)
\end{displaymath}
should be replaced by $\bar{\nu}_{in}$, and the $1+$ symbols in the terms
\begin{eqnarray*}
\left(1+\frac{\nu_{i\cdot}-\nu_{in}}{\bar{\nu}_{in}+r}\right) & \mathrm{and} & \left(1+\frac{\nu_{\cdot n}-\nu_{in}}{\bar{\nu}_{in}+r}\right)
\end{eqnarray*}
should be replaced by:
\begin{eqnarray*}
\left(1-\frac{\nu_{i\cdot}-\nu_{in}}{\bar{\nu}_{in}+r}\right) & \mathrm{and} & \left(1-\frac{\nu_{\cdot n}-\nu_{in}}{\bar{\nu}_{in}+r}\right)
\end{eqnarray*}
\end{enumerate}

\section*{\noindent Appendix B---Miscellaneous Proofs }

$ $

{\noindent \textbf{Proof of Lemma~\ref{lemma:1}:}}
To begin, recall that since $Z = X \times Y$, we can write $\rho_Z$ as a matrix
of real numbers, $\{\rho(x, y) : x \in X, y \in Y\}$ or, alternatively, as $(\rho_X, \rho_{Y \mid X})$,
where $\rho_{Y \mid X}$ is the set of $|X||Y|$ real numbers given by $\rho(x, y) / \rho_X(x)$ for all
$x \in X, y \in Y$. In other words, the space of all pairs $(\rho_{X}, \rho_{Y \mid X})$ is a coordinate system of $\Delta_{X \times Y}$,
given by the $|X| - 1$ real numbers specifying $\rho_X$ and the $|X|(|Y| - 1)$ real numbers specifying $\rho_{Y \mid X}$, as is the matrix
given by $|X||Y| - 1$ real numbers. Writing the coordinate transformation between
these coordinate systems as $\rho_{X,Y}(x, y) = \rho_{X}(x) \rho_{Y \mid X}(y \mid x)$,
we see that there is an integrating factor, which we can write as:
\begin{eqnarray}
\int d\rho_{X,Y} \ldots &=& \int d\rho_{X,Y} d\rho_X \prod_x \rho_X(x)^{|Y|-1} \ldots \nonumber
\end{eqnarray}
where we subtract one from the exponent inside the product to reflect the normalization
constraint on~$\rho_{Y \mid X}$.

The proposition's hypothesized equality expands to:
\begin{eqnarray}
\int d\rho_X \; Q(\rho_X) \frac{\prod_x \rho_X(x)^{\n_X(x) -1 + c / |X|}} {C(c, \n_X)}
 &=& 
 \int d\rho_{X,Y} \; Q(\rho_X) \frac{\prod_{x,y} \rho_{X,Y}(x,y)^{\n(x,y) -1 + c' / |X||Y|}} {C(c', \n_{X,Y})} \nonumber 
\end{eqnarray}
Using the appropriate integrating factor, the RHS can be written as:
\begin{eqnarray}
\!\!\!\!\!\! \!\!\!\!\!\! \!\!\!\!\!\! \!\!\!\!\!\! 
\int d\rho_X d\rho_{Y \mid X} \; Q(\rho_X) \prod_x \rho_X(x)^{|Y| - 1} \frac{\prod_{x,y} [\rho_X(x) \rho_{Y|X}(y \mid x)]^{\n(x,y) -1 + c' / |X||Y|}} {C(c', \n_{X,Y})} && \nonumber
\end{eqnarray}
{\small
\begin{eqnarray}
 \;\;\;\;\; &=& \int d\rho_X d\rho_{Y \mid X} \; Q(\rho_X) \prod_x \rho_X(x)^{|Y| - 1} \frac{\prod_{x,y} [\rho_X(x)]^{\n(x,y) -1 + c' / |X||Y|} \prod_{x,y} [\rho_{Y|X}(y \mid x)]^{\n(x,y) -1 + c' / |X||Y|}} {C(c', \n_{X,Y})} \nonumber \\
\end{eqnarray}}

Rearranging the exponents in the products and then separately collecting all terms involving $\rho_X$ and all
terms involving $\rho_{Y \mid X}$, we can rewrite this as
\begin{eqnarray}
&=& \int d\rho_X d\rho_{Y \mid X} \; Q(\rho_X) \frac{\prod_{x} [\rho_X(x)]^{n_X(x) - 1 + c' / |X|} \prod_{x,y} [\rho_{Y|X}(y \mid x)]^{\n(x,y) -1 + c' / |X||Y|}}
 {C(c', \n_{X,Y})} \nonumber \\
 &=& \int d\rho_X \frac{Q(\rho_X)}{C(c', \n_{X,Y})} \prod_{x} [\rho_X(x)]^{n_X(x) - 1 + c' / |X|} \prod_{x} \int d\rho_{Y \mid X} \; \prod_y [\rho_{Y|X}(y \mid x)]^{\n(x,y) -1 + c' / |X||Y|}
 \nonumber 
\end{eqnarray}
The inner integral is over the $|Y|$-dimensional simplex and evaluates to $\frac{\prod_y \Gamma(\n(x,y) + 
 c' / |X|Y|)} {\Gamma(\n_X(x) + |Y|)}$.
Therefore, rearranging terms, we get:
\begin{eqnarray}
\frac{1}{C(c', \n_{X,Y})} \frac{\prod_{x,y} \Gamma(\n(x,y) + c' / |X|Y|)}{\prod_x \Gamma(\n_X(x) + |Y|)}
 \int d\rho_X \; Q(\rho_X) \prod_{x} [\rho_X(x)]^{n_X(x) - 1 + c' / |X|} \nonumber
\end{eqnarray}
\begin{eqnarray}
 \;\;\;\;\; &=& \frac{\Gamma(N + c')}{\prod_x \Gamma(\n_X(x) + |Y|)}
 \int d\rho_X \; Q(\rho_X) \prod_{x} [\rho_X(x)]^{n_X(x) - 1| + c' / |X|} \nonumber
\end{eqnarray}
\begin{eqnarray}
 \;\;\;\;\; \;\;\;\;\; \;\;\;\;\; \;\;\;\;\; \;\;\;\;\; &=& \frac{1}{C(c', \n_X)} \int d\rho_X \; Q(\rho_X) \prod_{x} [\rho_X(x)]^{n_X(x) -1+ c' / |X|} \nonumber
\end{eqnarray}

By inspection, this expression equals the LHS of our hypothesized equality if $c = c'$. 
Going the other way, it is easy to see that if $Q(\rho_X) = \sum_x [\rho_X(x)]^2$, but $c \ne c'$, then our equality does not hold.
{\qed}

$ $

\noindent \textbf{Proof of Proposition~\ref{prop:iuv-weak}}
We have just proven that a DI hyperprior implies that IUV holds. To go the other way,
first, for any $\kk \in \Delta_X$, define:
\begin{eqnarray}
Q_{\kk}(\rho_X) &\equiv& \frac{ \delta(\rho_X - \kk) } { \prod_x k(x)^{n_X(x) - 1}} \nonumber
\end{eqnarray}

Next, plug Corollary~\ref{coroll:partial-iuv} into the definition of IUV to show that if IUV holds, then for any $Q$:
\begin{eqnarray}
\int dc \; F(c) \int d\rho_X \; Q(\rho_X) \D_{c, X}(\rho_X \; | \; \n_X) &=& 0 \nonumber
\end{eqnarray}
where $F(c)$ is defined as the analytic extension of $\pi(c \mid |X|) - \pi(c \mid |X||Y|)$.
Therefore, in particular, for any $\kk, c > 0$:
\begin{eqnarray}
0 &=& \int dc \; F(c) \int d\rho_X \; Q_{c, \kk}(\rho_X) \D_{c, X}(\rho_X \; | \; \n_X) \nonumber \\
&=& \int dc \frac{F(c) \bigg[\prod_x k(x)\bigg]^{c / |X|}} {C(c, \n_X)} \nonumber
\end{eqnarray}
Since $\prod_x k(x) \in [0, (1/|X|)^{|X|}]$, we see that for any $\alpha \in [0, 1/|X|]$:
 \begin{eqnarray}
\int dc \; B(c) \alpha^c &=& 0 \nonumber
\end{eqnarray}
where $B(c) \equiv \frac{F(c)} { C(c, \n_X)}$. Redefining $\alpha$ by dividing it by $(1 + \epsilon) |X|$ and
redefining $B(c)$ by multiplying it by $((1 + \epsilon) |X|)^c$, we see that for any $\alpha \in [0, (1 + \epsilon)]$:
 \begin{eqnarray}
\int dc \; B(c) \alpha^c &=& 0 \nonumber
\end{eqnarray}

Since, by hypothesis, this rescaled $B(c)$ is analytic about $c = 1$, we can differentiate both sides of this 
equation with respect to $\alpha$ an arbitrary number of times and evaluate
it at $\alpha = 1$. This establishes that all moments of $B(c)$ must equal zero. Since the Fourier
transform of $B$ is assumed analytic, this means that the Fourier transform of
$B(c)$ must equal zero identically. Therefore, $B(c)$ must equal zero identically, and~therefore, $F(c)$ must.

This establishes that if IUV holds then for any spaces, $X and Y$, it must be that $P(c \mid |X|) = P(c \mid |X||Y|)$. Relabeling
$X$ and $Y$ then establishes that if IUV holds, for any spaces, $X and Y$, $P(c \mid |X|) = P(c \mid |Y|)$.
{\qed}

$ $

{\noindent \textbf{Proof of Proposition~\ref{prop:new_like}:}}
Recall that $I(\SSS(\n) \subseteq Z)$ equals one iff $\SSS(\n)$, the support of $\n$ is a a subset
of $Z$. Therefore,
\begin{eqnarray}
P(\n \mid Z, m, \hat{Z}) &=& I(\SSS(\n) \subseteq Z) \int d\rho_Z \; P(\n \mid \rho_Z, Z, m, \hat{Z}) P(\rho_Z \mid Z, m, \hat{Z}) \nonumber \\
 &\propto& I(\SSS(\n) \subseteq Z) \int d\rho_Z \; \prod_{z \in Z} [\rho_Z(z)]^{n(z)} \D_{c,Z}(\rho_Z) \nonumber
\end{eqnarray}
Since we are using a Dirichlet prior with a uniform baseline distribution, by symmetry, 
the integral on the RHS must have the same value for all $Z$, such that $I(\SSS(\n)) \subseteq Z$ and $|Z| = m$. That value
is $G(\n, c, m)$.
This establishes the first claim. 

Next, write:
\begin{eqnarray}
P(\n \mid m, \hat{Z}) &=& \sum_{Z \subseteq \hat{Z}} P(\n \mid Z, m, \hat{Z})P(Z \mid m, \hat{Z}) . \nonumber
\end{eqnarray}
Combining this with our first result and with Equation~\eqref{eq:subset_sel_model}:
\begin{eqnarray}
P(\n \mid m, \hat{Z}) &\propto& \sum_{Z: |Z| = m} I(\SSS(\n) \subseteq Z) G(\n, c, m) \nonumber
\end{eqnarray}
(The fact that $P(Z \mid m, \hat{Z})$ is uniform over all $Z$ of size $m$ means that it will cancel out once we
divide by the appropriate sum to normalize $P(\n \mid m, \hat{Z})$.) 
For any set, $\SSS(\n)$, and any integer, $m \in \{|\SSS(\n)|, \ldots |\hat{Z}|\}$, there are a total of $\binom{|\hat{Z}| - |\SSS(\n)|}{m - |\SSS(\n)|}$ 
sets $Z \in \hat{Z}$, such that $I(\SSS(\n)) \subseteq Z$ and $|Z| = m$. Combining establishes the second claim.
{\qed}


\begin{thebibliography}{-------}
\bibitem{davidone}
Note that no particular property is required of the relation among multiple unseen random variables that may (or may not) exist.
\bibitem{davidtwo}
The proof of this proposition uses moment-generating functions with Fourier decompositions of the prior, $\pi$. To ensure we do not divide by zero, we have to introduce the constant $1 + \epsilon$ into that proof, and to ensure the convergence of our resultant Taylor decompositions, we have to assume infinite differentiability. This is the reason for the peculiar technical condition.
\bibitem{davidthree}
This is an example of a more general phenomenon, that in some statistical scenarios there is a low probability that the likelihood of a randomly sampled dataset is concentrated near the truth (various ``paradoxes'' of statistics can be traced to this phenomenon). See also the discussion of Bayes factors below.
\bibitem{davidfour}
Intuitively, for that $P(c)$, the likelihood says that a dataset, \{125, 125,\ldots \}, would imply a relatively high value of $c$ and, therefore, a high probability that $\rho$ is close to uniform over all bins. Given this, it also implies a low value of $|Z|$,  since if there were any more bins than the eight that have counts, we almost definitely would have seen them seen them for almost-uniform $\rho$. Conversely,  for the dataset, \{691, \ldots, 6\}, the implication is that $c$ must be low where some rare bins trail out, and as a result, there might be a few more bins who had zero counts.
\bibitem{davidfive}
Of course, since the formulas in WW implicitly assume $c = m$, care must be taken to insert appropriate pseudo-counts into those formulas if we want to use a value, $c$, that differs from $m$.
\end{thebibliography}
\end{document}